\documentclass[a4paper,12pt]{article}
\usepackage{dsfont}
\usepackage{amsmath,comment}

\usepackage{amssymb,amsthm}
\usepackage{mathrsfs,amsfonts}
\usepackage{color}
\usepackage{graphicx}
\usepackage{bbm}
\usepackage{ulem}
\usepackage[hypertexnames=false]{hyperref}

\allowdisplaybreaks

\allowdisplaybreaks

\newtheorem{con1}{Condition}[section]
\newtheorem{thm1}{Theorem}[section]
\newtheorem{def1}[thm1]{Definition}
\newtheorem{lem1}[thm1]{Lemma}
\newtheorem{cor1}[thm1]{Corollary}
\newtheorem{pro1}[thm1]{Proposition}
\newtheorem{rem1}[thm1]{Remark}
\newtheorem{ass1}[thm1]{Assumption}
\newtheorem{exa1}{\it{Example}}[section]

\numberwithin{equation}{section}

\def\bglemma{\begin{lem1}}\def\edlemma{\end{lem1}}

\def\bgproposition{\begin{pro1}}\def\edproposition{\end{pro1}}

\def\benumerate{\begin{enumerate}}\def\eenumerate{\end{enumerate}}
\def\bitemize{\begin{itemize}}\def\eitemize{\end{itemize}}

\def\beqlb{\begin{eqnarray}}\def\eeqlb{\end{eqnarray}}
\def\beqnn{\begin{eqnarray*}}\def\eeqnn{\end{eqnarray*}}

\DeclareMathOperator{\sgn}{sgn}\DeclareMathOperator{\var}{var}

\def\eqref#1{{\rm(\ref{#1})}}

\def\ar{\!\!\!&}

\def\proof{\noindent{\it Proof.~}}
\def\qed{\hfill$\square$\smallskip}

\def\mrm{\mathrm}\def\mbb{\mathbf}
\def\mbb{\mathbb}\def\mds{\mathds}

\def\d{\mrm{d}}\def\e{\mrm{e}}

\def\I{\mds{1}}

\newcommand{\R}{\mathbb R}
\newcommand{\pP}{\mathbb P}
\newcommand{\E}{\mathbb E}
\newcommand{\N}{\mathbb N}
\newcommand{\D}{\mathbb D}

\begin{document}

\title{Mean-field branching SDEs: propagation of chaos, scaling limits and phase transitions}

\author{Shukai Chen, Lina Ji, Xiaowen Zhou}
\date{}
\maketitle

\begin{abstract}
We study branching SDEs with law-dependent immigration and their mean-field particle approximations. Under a dissipativity condition and sufficiently weak interaction, a uniform propagation-of-chaos bound in time of order $N^{-1/2}$ is established.  On every fixed finite time horizon, the same order of propagation of chaos holds for arbitrary finite interaction strength. A two-stage scaling limit connects continuous-time discrete-state mean-field birth--death processes to interacting branching diffusions and then to the nonlinear equation.   For a logistic mean-field diffusion we prove a sharp criterion for extinction/non-extinction, and further show that  weak enough interaction strength is necessary for a uniform-in-time approximation.
\end{abstract}

\noindent\textbf{AMS 2020 Math Subject Classification.}
	60J80; 60K35; 60F05.

	\medskip
	

\section{Introduction and main results}\label{sec:introduction}

\subsection{Background and motivation}

Classical branching processes provide a fundamental probabilistic framework
for population dynamics.  Under high-density scaling, discrete branching models, such as Galton--Watson and continuous-time birth--death processes, converge to continuous-state branching diffusions; see, e.g., Feller \cite{Fel51} and Pardoux \cite{Pardoux16}. The classical example is the {\it Feller branching diffusion} introduced
in~\cite{Fel51}, which can be represented as
\begin{equation}\label{Feller}
X_t=X_0-b\int_0^tX_s \mathrm ds
+\int_0^t\int_0^{X_s}W(\mathrm ds,\mathrm du),
\end{equation}
where $W$ is space-time white noise. For models with immigration, Kawazu and
Watanabe~\cite{KW71} and Li~\cite{L06} established scaling limits to
continuous-state branching processes with immigration. The long-time behavior of branching diffusions is also well understood.
The Feller branching diffusion becomes extinct
almost surely when $b > 0$, whereas with nontrivial immigration it
admits a nondegenerate stationary distribution; see, e.g., Li \cite{L20}. For further background on Galton--Watson processes and continuous-state
branching processes, readers may refer to Athreya and Ney~\cite{AN72},
Li~\cite{L11}, the monograph~\cite{Pardoux16}, as well as the references therein.

The classical branching property relies on the independence of different
families and therefore does not capture interactions caused by limited
resources. After suitable scaling of discrete population models, branching
diffusions with nonlinear population regulation naturally arise in the form
\beqnn
X_t=x-b\int_0^tX_s\,\mathrm ds-\int_0^t g(X_s)\,\mathrm ds
+\int_0^t\int_0^{X_s}W(\mathrm ds,\mathrm du),
\eeqnn
see, e.g., \cite{Pardoux16}. Here, $g$ describes a population-dependent
interaction mechanism, such as competition. The logistic branching process was studied by Lambert~\cite{Lambert05}, while more general continuous-state
branching models with competition were investigated by Berestycki et al.~\cite{BFF18}. Furthermore, Li et al.~\cite{LYZ19a} studied a general class of
continuous-state nonlinear branching processes with population-size-dependent
branching rates and competition; Li et al.~\cite{LYZ19} also provide a discrete approximation. The long-time behavior of such models has also been studied. Extinction criteria were obtained in~\cite{Lambert05} and~\cite{LYZ19a}, while, in the presence of immigration, Li et al.~\cite{LLWZ25} proved exponential convergence to a unique stationary distribution under suitable conditions; see also Chen et al. \cite{CFJW26}.

 A different type of interaction arises in mean-field particle systems, where
each component depends on the empirical distribution of the entire system.
In the many-particle limit, this empirical distribution is replaced by the law
of a representative component, leading to McKean--Vlasov or
distribution-dependent stochastic equations. The notion of chaos originates
with Kac~\cite{Kac56}, while the probabilistic formulation of nonlinear
mean-field dynamics was introduced by McKean~\cite{Mckean66}. The McKean--Vlasov
limits and propagation of chaos were subsequently developed by
G\"{a}rtner~\cite{Gartner88} and Sznitman~\cite{Sznitman1991}. Related
distribution-dependent equations with square-root or irregular coefficients
have been studied by Bai and Huang~\cite{BH23} and Bao and
Huang~\cite{BH22}.  Liang et al.~\cite{LMW21} obtained a uniform-in-time estimate in $W_1$ for the $k$-particle marginals under partial dissipativity and a nondegenerate additive noise component, and this setting excludes square-root branching degeneracy at the origin. See also Huang \cite{H24} for a diffusion system driven by interacting noise. 

  For particle-level branching models, Claisse et al.~\cite{CKT24} and Cao et al.~\cite{CRT26} established qualitative and quantitative weak propagation of chaos, respectively. Fontbona and
Mu\~noz-Hern\'andez~\cite{FMH22} proved a quantitative mean-field limit for
logistic binary branching. Related mean-field systems of interacting
branching diffusions were studied by Hutzenthaler and
Wakolbinger~\cite{HW07} and Hutzenthaler~\cite{Hutz12}, with particular
attention to their long-time behavior and finite-particle approximations.   

Motivated by these works, we study a continuous-state branching
population with law-dependent immigration. 
Let
$f:\mathbb R_+\to\mathbb R_+$ be Lipschitz,
$g:\mathbb R_+\to\mathbb R_+$ be continuous and non-decreasing with
$g(0)=0$, and take $b\in\mathbb R$ and $a,\sigma>0$. On a filtered probability
space $(\Omega, \mathcal{F}, (\mathcal{F}_t)_{t \ge 0}, \mathbb{P})$ satisfying the usual conditions, let $W(\mathrm \d s,\mathrm \ du)$ be an
$(\mathcal{F}_t)$-adapted Gaussian white noise with intensity
$2\sigma\,\d s\,\d u$. We consider the following
mean-field branching stochastic equation:
\begin{equation}\label{X}
\begin{aligned}
X_t=X_0
&-b\int_0^tX_s\,\mathrm ds+\int_0^t\int_0^{X_s}W(\mathrm ds,\mathrm du)
-\int_0^tg(X_s)\,\mathrm ds\\
&
+a\int_0^t\mathbb E[f(X_s)]\,\mathrm ds .
\end{aligned}
\end{equation}
Here the initial value $X_0$ has law $\mu_0$ with finite second moment.
The last term in \eqref{X} represents a law-dependent immigration mechanism,
which introduces a mean-field interaction into the branching population
dynamics. We adopt the standard notion of solutions to \eqref{X}; see, for example, Wang \cite[Definition 1.1]{Wang18}. A c\`adl\`ag adapted process $(X_t)_{t\ge0}$ with values in $\R_+$ is called a strong solution if
$$
\E\left[\int_0^t|g(X_s)+X_s+f(X_s)|\,\d s\right]<\infty
$$
for every $t\ge0$ and \eqref{X} holds $\pP$-a.s. Weak solutions are defined analogously on an extended probability space. The following standard well-posedness result will be used.

\begin{pro1}\label{exist and unique}
There is a unique nonnegative strong solution to \eqref{X} and it is weakly well posed for the prescribed initial law. Moreover, for every $T<\infty$, $\sup_{0\leq t\leq T}\E [X_t^2]<\infty$.
\end{pro1}

The proof is given in Section \ref{sec:wellposedness-proofs}.

The natural finite-population counterpart replaces
$\mathbb{E}[f(X_t)]$ by the empirical mean
\[
N^{-1}\sum_{j=1}^N f(X_t^{j,N}).
\]
Formally, as $N$ tends to infinity, this empirical quantity approaches
$\mathbb{E}[f(X_t)]$, and every fixed collection of particles becomes
asymptotically independent with common law
$\operatorname{Law}(X_t)$. This is precisely the propagation-of-chaos
principle studied below.

There is also a microscopic question behind the diffusion system.
Branching diffusions are themselves continuous-state approximations of
discrete population processes. Accordingly, the models in this paper
form the three-level hierarchy
\[
Z^{n,N}/n
\xrightarrow{\,n\to\infty,\;N\ \mathrm{fixed}\,}
X^N
\xrightarrow{\,N\to\infty\,}
X,
\]
where $Z^{n,N}/n$ is the rescaled discrete mean-field birth--death
population, $X^N$ is the interacting branching diffusion system, and
$X$ solves~\eqref{X}. Establishing both transitions gives a direct
connection between a microscopic discrete population model and the
nonlinear mean-field branching equation.

A particularly informative  case is the logistic choice
$f(x)=x$ and $g(x)=b_1x+b_2x^2$. Its phase transition distinguishes
finite-time mean-field approximation from approximation uniformly over an
infinite time horizon, and thereby tests whether the small-interaction
condition in our main theorem reflects the dynamics or only the proof method.

The main results answering these questions are
stated in the next subsection.

\subsection{Main results}

\subsubsection{Propagation of chaos}\label{subsec:particle-system}

We focus on the mean-field interacting particle system corresponding to \eqref{X} which is defined, for all $N\in\N=\{1,2,...\}$,
\beqlb\label{Xi}
X_t^{i,N} = X_0^{i,N}\ar+\ar
\int_0^t\int_0^{X_{s-}^{i,N}}W_i(\d s,\d u)-\int_0^t[bX_s^{i,N} + g(X_s^{i,N})]\,\d s+a\int_0^t \bar{f}^N_s\,\d s
\eeqlb
for each $1\le i\le N$, where
$$
\bar{f}^N_t:=\frac{1}{N}\sum_{j=1}^N f(X_t^{j,N}),
$$
$\{W_i(\d s,\d u): s>0, u>0\}_{1\le i\le N}$ are independent copies of $W$, and $X^{i,N}_0$ each have law $\mu^N_0$ with finite second moment. The generator of $\{(X^{i,N}_t)_{t\ge0}: 1\le i\le N\}$ is given by
\beqlb\label{generator of particle}
{\mathcal{L}_{N}}F(x):=\sum_{i=1}^N{\mathcal{L}^i_{N}}F(x),
\eeqlb
where
\beqlb\label{coupling 1.5}
{\mathcal{L}^i_{N}}F(x)
=\sigma x^iF''_{x^ix^i}(x)
+[-bx^i-g(x^i)]
F'_{x^i}(x)+\frac{a}{N}\sum_{j=1}^Nf(x^j)F'_{x^i}(x)
\eeqlb
for every twice continuously differentiable function $F:\R^N_+\mapsto\R_+$.

\begin{pro1}\label{prop:particle-wellposedness}
There is a pathwise unique nonnegative strong solution to the particle system \eqref{Xi} and it is
weakly well posed.  Moreover, for every $T<\infty$,
\[
\sup_{0\leq t\leq T}\E\left[\left(\sum_{i=1}^N X_t^{i,N}\right)^2\right]<\infty.
\]
\end{pro1}

The proof is given in Section~\ref{sec:wellposedness-proofs}.

The connection between \eqref{X} and \eqref{Xi} is that the 
mean-field-type SDE \eqref{X} describes the dynamics of one particle of the interacting particle system \eqref{Xi} when the number of particles $N$ tends to infinity.

Denote by $(X_t)_{t\ge0}$ and $\{(X^{i,N}_t)_{t\ge0}:1\le i\le N\}$ the unique strong solutions to \eqref{X} and \eqref{Xi} with initial distributions $\mu_0$ and $(\mu_0^N)^{\bigotimes N}$, respectively.

For any $d\ge1$, let $C^2(\R_+^d)$ be the linear space of twice continuously differentiable functions on $\R_+^d$. For any $h\in C^2(\R_+)$, we write
\beqlb\label{operator}
{\mathcal{L}}h(x)=\sigma xh''(x)
+[-bx-g(x)]h'(x).
\eeqlb

For any $N\ge1$, let ${\mathcal{P}}$ be the collection of all probability measures on $\R_+^N$ equipped with the weak topology. Let
$$
{\mathcal{P}}_{d^N}:=
\left\{\mu\in{\mathcal{P}}: \int_{\R_+^N}
d^N(x,0)\,\mu(\d x)<\infty
\right\},
$$
where the metric
\beqlb\label{metric}
d^N(x,y):=\frac{1}{N}
\sum_{i=1}^N|x^i-y^i|,\quad
x=(x^1,...,x^N), \quad y=(y^1,...,y^N)\in\R_+^N
\eeqlb
and $\R_+^N=\left\{x|x=(x^1,...,x^N), x^i\in\R_+\right\}$.
The associated {\it Wasserstein distance} is given by
$$
W_{d^N}(\gamma,\eta)
=\inf_{\pi\in{{\mathcal{C}}(\gamma,\eta)}}\int_{\R_+^{2N}}d^{N}(x,y)\,\pi(\d x,\d y),\quad \gamma,\eta\in{\mathcal{P}}_{d^N},
$$
where ${\mathcal{C}}(\gamma,\eta)$
is the collection of all probability measures on $\R_+^{2N}$ with marginals $\gamma$ and $\eta$.

Let $C_f>0$ be a Lipschitz constant of $f$.  For any $l>0$, define
\[
\eta(l):=\inf_{\substack{y\geq0\\ r\geq l}}
\left[b+\frac{g(y+r)-g(y)}{r}\right],
\qquad
c_1:=\frac{2|b|+1}{\sigma},
\]
and
\begin{equation}\label{eq:a-star}
a_*
:=\frac1{C_f}
\sup_{\{l>0:\eta(l)>0\}}
\frac{\min\{1,\eta(l)\}}
{(1+\e^{c_1l})^2}.
\end{equation}

\begin{con1}\label{con2}
\beqnn
\liminf_{r\rightarrow\infty}\inf_{y\ge0}\left[\
b+\frac{g(y+r)-g(y)}{r}\right]>0.
\eeqnn
\end{con1}

The Condition~\ref{con2} implies that $\eta(l)>0$ for all sufficiently large
$l$, and hence $a_*>0$.

\begin{thm1}\label{main result}
Suppose that Condition \ref{con2} holds.  Then, for every
$0<A<a_*$, there exist strictly positive constants
$C_1,C_2$ and $\lambda$, depending on $A$ and $\mu_0$, such that
for all $t\geq0$, $N\geq2$ and $a\in(0,A]$,
\begin{equation}\label{unifom chaos}
W_{d^N}(\mu^{\bigotimes N}_t, \mu_t^N) \leq C_1 e^{-\lambda t}W_{d^N}(\mu^{\bigotimes N}_0, (\mu_0^N)^{\bigotimes N}) + C_2 N^{-1/2},
\end{equation}
where $\mu^{\bigotimes N}_t$ is the product law of $N$ independent solutions to \eqref{X} with initial distribution $\mu_0$, and $\mu^N_t$ is the law of the particles driven by \eqref{Xi} with  initial distribution $(\mu_0^N)^{\bigotimes N}$.
\end{thm1}

  The proof of Theorem \ref{main result} is given in Section~\ref{sec:main-proof}. 

\begin{def1}\label{def:thresholds}
We say that a constant $A>0$ is admissible if, for every initial law $\mu_0$
with finite second moment, there exist constants $C_1,C_2,\lambda>0$, allowed
to depend on $A$ and $\mu_0$ but not on $t$, $N$ or $a$, such that
\eqref{unifom chaos} holds for all $t\geq0$, $N\geq2$, $a\in(0,A]$ and every
initial law $\mu_0^N$ with finite second moment. Let
$\mathcal S_{\rm uc}$ be the set of all admissible $A$ and define
\begin{equation}\label{eq:a-critical-general}
a_c :=\sup\mathcal S_{\rm uc}\in[0,\infty].
\end{equation}
\end{def1}

The threshold $a_c$ is intrinsic. It records the supremal range of interaction
strengths on which a uniform-in-time propagation-of-chaos estimate of the
above form can hold. In contrast, the threshold $a_*$ is explicit and proof-dependent. It follows from Definition \ref{def:thresholds} that

\begin{cor1}\label{corollary}
Suppose that Condition \ref{con2} holds. Then we have
$a_*\le a_c$.
\end{cor1}

Let $b^-:=(-b)\vee0$ and
\[
\Phi_\kappa(t):=
\begin{cases}
(\e^{\kappa t}-1)/\kappa,& \kappa>0,\\
t,& \kappa=0.
\end{cases}
\]
Theorem \ref{main result} requires a smallness
condition on $a$.  The next theorem shows that, on a fixed finite time horizon, the same order of
propagation of chaos holds for every fixed interaction strength; this estimate
will also be used in the scaling limit below.
For any \(0<T<\infty\), we define
$$
M_{T,a}:=\sup_{0\leq s\leq T}(\E[f(X_s)^2])^{1/2},\quad \kappa_a:=b^-+aC_f.
$$
\begin{thm1}\label{propagation 2}
Fix $T<\infty$. Then for all $t\in[0,T]$ and $N\geq2$, we have 
\[
W_{d^N}(\mu^{\bigotimes N}_t,\mu_t^N)
\leq
\e^{\kappa_a t}W_{d^N}(\mu^{\bigotimes N}_0,(\mu_0^N)^{\bigotimes N})
+aM_{T,a}\Phi_{\kappa_a}(t)N^{-1/2}.
\]
\end{thm1}

 The proof of Theorem \ref{propagation 2} is deferred to Section   ~\ref{sec:main-proof}. 

\subsubsection{Scaling limit of a discrete mean-field population}\label{subsec:scaling-result}
 We next give a discrete population model. For $r\in\R$, write $r^+=\max\{r,0\}$ and $r^-=\max\{-r,0\}.$
Let $\mathbb Z_+=\{0,1,2,\ldots\}$. For $n,N\geq1$, let $Z^{n,N}=(Z^{1,n,N},\ldots,Z^{N,n,N})$
be a continuous-time Markov chain on $\mathbb Z_+^N:=\left\{z|z=(z^1,...,z^N),~z^i\in\mathbb Z_+\right\}$. For $z\in\mathbb Z_+^N$, define
\[
f^{n,N}(z):=\frac1N\sum_{j=1}^N f(z^j/n).
\]
  Fix $0<A<a_*$ and let $a\in(0,A]$. 
If $Z_t^{n,N}=z$, then its $i$th coordinate has the transitions
\begin{align}
z^i&\longrightarrow z^i+1
&&\text{at rate }\lambda_i^{+,n,N}(z)
=z^i(n\sigma+b^-)+na f^{n,N}(z),                                      \label{eq:discrete-birth-rate}\\
z^i&\longrightarrow z^i-1
&&\text{at rate }\lambda_i^{-,n,N}(z)
=z^i(n\sigma+b^+)+ng(z^i/n).                                          \label{eq:discrete-death-rate}
\end{align}
The nearly equal birth and death rates $z^i(n\sigma+b^\pm)$ give the diffusion term in \eqref{Xi}. The extra
death rate $ng(z^i/n)$ gives the competition term, and the extra birth rate $naf^{n,N}(z)$ gives the mean-field term.
Our second main result states that the rescaled chain converges first to the particle system
and then, as the number of particles tends to infinity, to \eqref{X}.

\begin{thm1}\label{thm:discrete-approximation}
Suppose that Condition \ref{con2} holds.
 Let $\xi_1,\xi_2,\ldots$ be i.i.d. with law $\mu_0$. Set $Z_0^{i,n,N}=\lfloor n\xi_i\rfloor$ and $X_0^{i,n,N}=\frac{Z_0^{i,n,N}}n.$
Let $Z^{n,N}$ be the chain defined by
\eqref{eq:discrete-birth-rate}--\eqref{eq:discrete-death-rate} and let
$X^{n,N}=Z^{n,N}/n$.  Then, for every fixed $T<\infty$, the following hold.
\begin{itemize}
\item[(i)] For each fixed $N$,
\[
X^{n,N}\Rightarrow X^N\qquad\text{in }\D([0,T],\R_+^N),\qquad n\to\infty,
\]
where $X^N=(X^{1,N},\ldots,X^{N,N})$ solves \eqref{Xi} with initial state
$(\xi_1,\ldots,\xi_N)$.
\item[(ii)] The iterated limit satisfies
\[
\lim_{N\to\infty}\lim_{n\to\infty}\operatorname{Law}(X^{1,n,N})
=\operatorname{Law}(X)
\]
weakly on $\D([0,T],\R_+)$, where $X$ solves \eqref{X} with initial law $\mu_0$.
\item[(iii)] There is a deterministic sequence $N_T(n)\to\infty$ such that
\[
X^{1,n,N_T(n)}\Rightarrow X\qquad\text{in }\D([0,T],\R_+),\qquad n\to\infty.
\]
\end{itemize}
\end{thm1}

The proof of Theorem \ref{thm:discrete-approximation} is deferred to Section \ref{sec:discrete-approximation}.

\subsubsection{A logistic mean-field branching diffusion}\label{subsec:logistic-model}

To compare the explicit sufficient threshold $a_*$ with the intrinsic
threshold $a_c$, we now specialize in a logistic model for which the
long-time behavior is known. Take $f(x)=x$ and
$g(x)=b_1x+b_2x^2$, where $b_1\geq0$ and $b_2>0$. Then \eqref{X} becomes
\beqlb\label{eq:logistic-mv}
X_t=X_0+\int_0^t\bigl[a\mathbb E[X_s]-bX_s-b_1X_s-b_2X_s^2\bigr]\,\d s
+\int_0^t\int_0^{X_s}W(\d s,\d u).
\eeqlb
In population terms, quadratic drift $-b_2X_s^2$ models local
density-dependent competition, while $a\mathbb E[X_s]$ is an input to immigration
or recruitment generated by the mean population size. When $a=0$,
\eqref{eq:logistic-mv} reduces to the continuous-state logistic branching 
process studied in~\cite{Lambert05}. Locally regulated branching
populations and their mean-field and island approximations were studied
in~\cite{HW07} and~\cite{Hutz12}; see also~\cite{FMH22}.

Set
$\beta=b+b_1$. For this choice of $f$ and $g$, the explicit threshold \eqref{eq:a-star}
becomes
\[
a_*
=\sup_{l>0}
\frac{\min\{1,(\beta+b_2l)^+\}}{(1+\e^{c_1l})^2},
\qquad c_1=\frac{2|b|+1}{\sigma},
\]
since $C_f=1$ and
$\eta(l)=\inf_{y\geq0,r\geq l}\{\beta+b_2(2y+r)\}=\beta+b_2l$. For $\nu>0$, we let
\beqlb\label{scale fun1}
U(x)=\frac{\beta}{\sigma}x+\frac{b_2}{2\sigma}x^2,\qquad w(x)=\e^{-U(x)},
\qquad
I_\nu=\int_0^\infty x^{\nu-1}w(x)\,\d x
\eeqlb
and define
\beqlb\label{eq:critical-ac}
a_\star^{\rm log}:=\frac{\sigma}{I_1}
=\sigma\left\{\int_0^\infty
\exp\left(-\frac{\beta}{\sigma}y-\frac{b_2}{2\sigma}y^2\right)\,\d y\right\}^{-1}.
\eeqlb
Recall the intrinsic constant \(a_c\) defined in \eqref{eq:a-critical-general}.  The proofs of the following results are deferred to Section \ref{example}.  

\begin{pro1}\label{prop: bound of ac}
For the mean-field logistic SDE \eqref{eq:logistic-mv}, we have
\[
a_*\le a_c\le a_\star^{\rm log},
\]
where $a_\star^{\rm log}$ is given by \eqref{eq:critical-ac}.
\end{pro1}

\begin{rem1}\label{rem:upper-bound-intuition}
The upper bound $a_c\leq a_\star^{\rm log}$ reflects a noncommutation between
the long-time and mean-field limits. When $a_c>a_\star^{\rm log}$, the nonlinear
equation admits a nonzero invariant law, as stated in
Proposition~\ref{prop:logistic-threshold} \textup {(ii)}. If initialized with
this law, the limiting equation remains stationary. In contrast, for every
fixed $N$, the corresponding particle system is eventually absorbed at the
zero configuration. Therefore, the long-time Wasserstein distance between the
two systems stays bounded away from zero, whereas uniform-in-time propagation
of chaos would require it to vanish as $N\to\infty$. Hence an admissible
interaction range cannot extend beyond $a_\star^{\rm log}$.
\end{rem1}

  The quantity $a_\star^{\rm log}$ is also the critical value for the long-term
behavior of the nonlinear population law, so we record the corresponding phase
transition for later use.

\begin{pro1}\phantomsection\label{prop:logistic-threshold}
\begin{itemize}
\item[(i)] If $a\leq a_\star^{\rm log}$, then $X_t\to0$ in $L^1$.
\item[(ii)] If $a>a_\star^{\rm log}$, there is a unique $\alpha_a>0$ satisfying $a=\sigma\alpha_a I_{\alpha_a}/I_{\alpha_a+1}$. The probability measure
$\pi_a(\d x)=\frac{x^{\alpha_a-1}w(x)}{I_{\alpha_a}}\d x$
has moments of all orders and is the unique nonzero invariant probability
measure of \eqref{eq:logistic-mv} with finite mean.
\end{itemize}
\end{pro1}

\begin{rem1}\label{rem:hw-logistic-threshold}
The phase-transition and invariant-law classification in
Proposition~\ref{prop:logistic-threshold} are essentially contained
in~\cite[Lemmas~5.1--5.2 and Corollary~6]{HW07}. Indeed, under
some proper parameter transformation,
the mean-field equation~\cite[(69)]{HW07} coincides with
\eqref{eq:logistic-mv}, and their criticality condition~\cite[(84)]{HW07}
becomes $aI_1/\sigma=1$, equivalently $a=a_\star^{\rm log}$. Their extinction conclusion is weak convergence to $\delta_0$; together with
the uniform second-moment bound in
Lemma~\ref{lem:logistic-moments}, this yields the $L^1$ convergence in part~\textup{(i)}.  We restate and prove the result in our notation to make the
parameter correspondence explicit and because both parts are used below:
part~\textup{(i)} in the subcritical boundary analysis and
part~\textup{(ii)} in the proof of $a_c\leq a_\star^{\rm log}$. Thus the
proposition is included as an auxiliary result rather than claimed as a new
contribution. 
\end{rem1}

Proposition~\ref{prop:logistic-threshold} identifies the nonzero invariant law
in the supercritical regime, but it does not by itself give convergence to
this law from a deterministic positive initial state. This convergence,
together with the convergence of the first moment, is needed in the following because the
asymptotic mean-field immigration rate $am_t$ enters the boundary
classification. The following result presents these long-time limits.


\begin{pro1}\label{p0814}
Assume that $a>a_\star^{\rm log}$ and let $m_t = \mbb{E}_x[X_t]$. Then, for every $x>0$, as $t \rightarrow \infty$, we have $X_t\ \Rightarrow\ \pi_a$ and
\beqnn
m_t\longrightarrow m_a
=\int_0^\infty y\,\pi_a(\d y) =\frac{I_{\alpha_a+1}}{I_{\alpha_a}}
=\frac{\sigma\alpha_a}{a}.
\eeqnn
\end{pro1}

The invariant-law classification does not determine whether the population started from a positive size can reach the extinction boundary before
mean-field immigration restores positive mass.
For a deterministic initial value $x>0$, write $\mbb P_x$ and $\mbb E_x$
for the corresponding probability and expectation, and let
$m_t:=\mbb E_x[X_t]$ and $\tau_0:=\inf\{t\geq0:X_t=0\}$, and define $a_{\text{hit}}:=\frac{\sigma I_1}{I_2} $.
The next theorem describes this boundary behavior.

\begin{thm1}\label{thm:logistic-boundary}
Let $X$ be the solution of the logistic mean-field equation displayed above
with $X_0=x>0$. Then the following assertions hold.
\begin{itemize}
\item[(i)] (Exact no-hitting criterion) $\mbb P_x(\tau_0<\infty)=0$ if and only if $am_t\geq\sigma$ for every $t\geq0$.
\item[(ii)] (Almost-sure hitting under eventual subcritical immigration) If $\limsup_{t\to\infty}am_t<\sigma$, then
$\mbb P_x(\tau_0<\infty)=1$. In particular, this conclusion holds when $a\leq a_\star^{\rm log}$;
\item[(iii)] (Classification by the interaction strength) $ a_{\text{hit}}>a_\star^{\rm log}$. Moreover,
$\mbb P_x(\tau_0<\infty)=1$ if $a<a_{\text{hit}}$, whereas
$\mbb P_x(\tau_0<\infty)<1$ if $a>a_{\text{hit}}$.
\end{itemize}
\end{thm1}

\begin{rem1}\label{rem:0}
Combining ~(i) and (iii) of Theorem \ref{thm:logistic-boundary}, if $a>a_{\text{hit}}$ and $am_t<\sigma$ for some $t\geq0$,
then $0<\mbb P_x(\tau_0<\infty)<1$. In particular, this holds whenever
$a>a_{\text{hit}}$ and $0<x<\sigma/a$.
\end{rem1}

\begin{rem1}\label{rem:logistic-boundary-interpretation}

Theorem \ref{thm:logistic-boundary}~\textup{(ii)} and~\textup{(iii)} are linked through the long-time behavior of the deterministic mean. More precisely, when $a>a_\star^{\rm log}$, Proposition~\ref{p0814} gives $am_t\longrightarrow \sigma\alpha_a$ for $F(\alpha_a)=a$. Since $F$ is strictly increasing and $F(1)=a_{\text{hit}}$, we have
\[
\alpha_a<1 \ \Longleftrightarrow\ a<a_{\text{hit}},
\qquad
\alpha_a>1 \ \Longleftrightarrow\ a>a_{\text{hit}}.
\]
Consequently, for $a_\star^{\rm log}<a<a_{\text{hit}}$ we have $\limsup_{t\to\infty}am_t<\sigma$ and Theorem\ref{thm:logistic-boundary}~\textup{(iii)} thus follows directly from Theorem\ref{thm:logistic-boundary}~\textup{(ii)}. For $a>a_{\text{hit}}$, the limiting mean immigration rate is instead larger than $\sigma$, so the hypothesis of Theorem\ref{thm:logistic-boundary}~\textup{(ii)} no longer applies; Theorem\ref{thm:logistic-boundary}~\textup{(iii)} then gives the complementary conclusion that boundary hitting is not almost sure. Together with the case $a\le a_\star^{\rm log}$ already covered by Theorem\ref{thm:logistic-boundary}~\textup{(ii)}, this identifies $a_{\text{hit}}$   as the boundary-hitting threshold away from the critical equality.
\end{rem1}

\begin{rem1}
We point out that the boundary point $0$ is not absorbing for \eqref{eq:logistic-mv}. Indeed, when $X_t=0$, the branching noise term $\int_0^{X_t}W(\d s,\d u)$ vanishes, whereas the drift still contains the immigration term $a\,m_t$. Consequently, as long as $m_t>0$, any trajectory reaching $0$ is pushed back into $\mathbb R_+$ by the positive mean-field immigration. Only the degenerate solution started from $\delta_0$, for which $m_t\equiv0$, remains identically equal to $0$.
\end{rem1}

\subsection{Contributions and proof strategy}\label{subsec:contributions}

The estimate in Theorem~\ref{main result} has two components. The
exponentially decaying term is the ergodic component, which describes the loss of
dependence on the initial distributions. The negative drift generated by
competition at large distances is favorable here because it provides the dissipativity required for uniform moment bounds and exponential contraction. The term of order $N^{-1/2}$ is the propagation of chaos error caused by
replacing the law-dependent mean with an empirical average.

We estimate this
term by coupling each particle with an independent copy of the limiting
process. The empirical discrepancy is decomposed into a centered fluctuation
and a coupling error, and the latter is absorbed by the contraction when $a$
is sufficiently small. A concave distance is used to handle the non-Lipschitz
square-root coefficient near the origin. Consequently, if
$\mu_0^N=\mu_0$, the approximation error is of the order $N^{-1/2}$ uniformly for all $t\geq0$.

Theorem~\ref{thm:discrete-approximation} connects a microscopic mean-field
birth--death process to \eqref{X}. Rescaling time or state is a standard way
to relate discrete and continuous branching models; see
\cite{KW71}, \cite{LYZ19}, \cite{L06} and the references therein. Here, the three assertions
distinguish the two limiting operations. Part~(i) gives the diffusion
approximation for fixed $N$, part~(ii) subsequently takes the mean-field limit,
and part~(iii) realizes both limits along a deterministic regime $N=N_T(n)$.
To the best of our knowledge, such a scaling limit has not previously been
established for a branching population model with the distribution-dependent
interaction considered here.

For each fixed $N$, the diffusion approximation is established by proving
tightness and convergence of the generators and then identifying the limiting
martingale problem. Propagation of chaos then shows that, as $N\to\infty$, the law of any
fixed component converges to that of the McKean--Vlasov solution. A diagonal argument then yields the sequence
$N_T(n)$. Since the fixed-$N$ diffusion approximation is qualitative and
provides no convergence rate in $n$, the argument does not give an explicit
relation between $N$ and $n$.

For the diffusion of the logistic mean-field branch, 
Theorem 1.10 gives a sharp boundary-hitting classification for the diffusion of the logistic mean-field branch.
It's proof combines Feller's boundary test for the diffusion with frozen immigration, the long-time limit of $m_t$, and the monotonicity of $F(\alpha)=\sigma\alpha I_\alpha/I_{\alpha+1}$.
Another  contribution is the upper bound
$a_c\leq a_\star^{\rm log}$ in Proposition~\ref{prop: bound of ac}. Together
with Theorem~\ref{main result}, it places the unknown intrinsic threshold
$a_c$ in the explicit interval $[a_*,a_\star^{\rm log}]$.
The proof uses Proposition~\ref{prop:logistic-threshold} (ii) by initializing
the limiting equation with its nonzero invariant law $\pi_a$ for
$a>a_\star^{\rm log}$. The limiting law then remains stationary. By contrast,
the total mass of every fixed-$N$ particle system is dominated by a logistic
branching process that is eventually absorbed at zero. This difference yields
a positive long-time Wasserstein separation and rules out a uniform estimate
with an error that disappears as $N\to\infty$.

\section{Proof of Theorem \ref{main result} and Theorem \ref{propagation 2}}\label{sec:main-proof}

\subsection{The decoupled equation}

For a fixed initial law, the curve $t\mapsto\E[f(X_t)]$ is deterministic, so
$(X_t)_{t\geq0}$ is a time-inhomogeneous Markov process.  What fails in general is a
homogeneous transition semigroup independent of the initial law.  For
the coupling argument, it is nevertheless convenient to make the deterministic law curve
explicit by freezing it.
We consider the decoupled stochastic equation associated with equation \eqref{X}:
\beqlb\label{decoupled SDE}
\bar{X}_t=\bar{X}_0+\int_0^t\int_0^{\bar{X}_{s-}}W(\d s,\d u)-\int_0^t[b\bar{X}_s+ g(\bar{X}_s)]\,
\d s +a \int_0^t\E[ f(X_s)]\,\d s
\eeqlb
with $Law(\bar{X}_0)=\mu_0$. Note that equation \eqref{decoupled SDE} is a classical stochastic equation (distribution-independent) with time-dependent drift. It is easy to see that  the stochastic equation \eqref{decoupled SDE} has a unique strong solution $(\bar{X}_t)_{t\ge0}$. The process $(\bar{X}_t)_{t\ge0}$ is a time-inhomogeneous Markov process and its time-dependent generator, acting on $h\in C^2(\R_+)$, is given by
$$
Q_{X_t}h(x)= {\mathcal{L}}h(x)
+a\E [f(X_t)]h'(x),
$$
where ${\mathcal{L}}$ is given by \eqref{operator}. Furthermore, since both \(X\) and \(\bar X\) solve \eqref{decoupled SDE} with the same initial condition and driving white noise, pathwise uniqueness implies that they are indistinguishable. Consequently, $\operatorname{Law}(\bar X_t)=\operatorname{Law}(X_t)=\mu_t$ for $t \ge 0$.
Thus, to prove Theorem~\ref{main result}, it remains to compare the particle system \eqref{Xi} with a system of independent copies of the decoupled process \eqref{decoupled SDE}.

\subsection{The coupling estimate}

Let $\pi_0$ be a coupling of $\mu_0$ and $\mu_0^N$ such that $\int_{\R_+^2} |x-y|\,\pi_0(\d x,\d y)=W_1(\mu_0,\mu_0^N).$
On a common probability space, let $\{(X_0^i,X_0^{i,N})\}_{1\le i\le N}$ be i.i.d. random variables with law $\pi_0$, independent of $\{W_i\}_{1\le i\le N}$, and let the particle system \eqref{Xi} start from $(X_0^{1,N},\ldots,X_0^{N,N})$. By the definition of $d^N$ in \eqref{metric}, we have
$$
\int_{\R_+^{2N}} d^N(x,y)\,\pi_0^{\otimes N}(\d x,\d y)
=W_{d^N}(\mu_0^{\otimes N},(\mu_0^N)^{\otimes N}).
$$
Hence, with $X_0=(X_0^1,\ldots,X_0^N)$ and $X_0^N=(X_0^{1,N},\ldots,X_0^{N,N})$,
\begin{equation}\label{eq:optimal-initial-coupling}
\E [d^N(X_0,X_0^N)]=W_{d^N}(\mu_0^{\otimes N},(\mu_0^N)^{\otimes N}).
\end{equation}

Let $\{(X^i_t)_{t\ge0}:1\le i\le N\}$ be $N$ independent versions of the solution to equation \eqref{X}, determined as
\beqnn
X^i_t = X^i_0\ar+\ar\int_0^t\int_0^{X^i_{s-}}W_i(\d s,\d u)-\int_0^t[bX_s^i + g(X^i_s)]\,
\d s +a \int_0^t\E [f(X^i_s)]\,\d s 
\eeqnn
for $1\le i\le N$,
where $X^i_0$ are i.i.d. random variables with $Law(X^i_0)=\mu_0$ for each $1\le i\le N$, and
$\{W_i(\d s, \d u):s>0, u>0\}_{1\le i\le N}$ are the same as those in \eqref{Xi}.
Then for each $1\le i\le N$, $(X^i_t)_{t\ge0}$ has the same law as the solution to equation \eqref{X}.
Furthermore, we consider the corresponding decoupled stochastic equation
\beqlb\label{coupling SDE}
\bar{X}^i_t\ar=\ar X^i_0+\int_0^t\int_0^{\bar{X}^i_{s-}}W_i(\d s,\d u)-\int_0^t[b\bar{X}^i_s + g(\bar{X}^i_s)]\,
\d s+a \int_0^t\E[f(X_s)]\,\d s 
\eeqlb
for $1 \le i \le N.$
Because the white noises $\{W_i(\d s,\d u):s>0, u>0\}$ are independent of each other, $(\bar{X}_t^i)_{t\ge0}$ are independent of each other, and $Law(\bar{X}^i_t)=\mu_t$ for each $1\le i\le N$ and $t>0$. The corresponding generator of $\{(\bar{X}^i_t)_{t\ge0}:1\le i\le N\}$ is given  by ${\mathcal{L}}_{X_t}H(x)=\sum_{i=1}^N{\mathcal{L}}^i_{X_t}H(x),$
where
\beqlb\label{coupling 2.11}
{\mathcal{L}}^i_{X_t}H(x)=\sigma x^i H''_{x^ix^i}(x)+[-bx^i-g(x^i)]H'_{x^i}(x)+a\E [f(X_t)] H'_{x^i}(x).
\eeqlb
Here, $H_{x}'$ denotes the first derivative with respect to $x$, $H_{xx}''$ denotes the second derivative with respect to $x$.
Therefore, in order to obtain Theorem \ref{main result}, it suffices to quantify the propagation of chaos of the  particle system  \eqref{Xi} towards the decoupled SDE \eqref{coupling SDE}.

For any $x=(x^1,...,x^N),\,y=(y^1,...,y^N)\in \R^{N}_+$ and the solution $X$ to \eqref{X}, the generator of the process $\{(\bar{X}^i_t, X^{i,N}_t)_{t\ge0}:1\le i\le N\}$, acting on the twice continuously differentiable function $F: \R^{2N}_+\mapsto \R$, is given by
\beqnn
{\bf L}_{X_t,N}F(x,y):=\sum_{i=1}^N{\mathcal{L}}^i_{X_t,N}F(x,y),
\eeqnn
where
\beqnn
{\mathcal{L}}^i_{X_t,N}F(x,y)=\tilde{\mathcal{L}}^iF(x,y)+a\E f(X_t)F'_{x^i}(x,y)+\frac{a}{N}\sum_{j=1}^Nf(y^j)F'_{y^i}(x,y)
\eeqnn
and
\beqlb\label{coupling 2.15}
\tilde{\mathcal{L}}^iF(x,y)\ar=\ar[-bx^i-g(x^i)]F'_{x^i}(x,y)+[-by^i-g(y^i)]F_{y^i}'(x,y)\cr
\ar\ar+\sigma x^iF''_{x^ix^i}(x,y)+\sigma y^iF_{y^iy^i}''(x,y)+2\sigma(x^i\wedge y^i)F''_{x^iy^i}(x,y).
\eeqlb
Clearly, ${\mathcal{L}}^i_{X_t,N}$ is the coupling generator of ${\mathcal{L}}^i_{X_t}$ and ${\mathcal{L}}^i_N$ given by \eqref{coupling 2.11} and \eqref{coupling 1.5} in the sense that
$
{\mathcal{L}}^i_{X_t,N}F_0(x,y)={\mathcal{L}}^i_{X_t}F_1(x)+{\mathcal{L}}^i_N F_2(y)$
for any $F_0(x,y)=F_1(x)+F_2(y)$.

The condition \ref{con2} implies that there exist constants $l_0 > 0$ and $\eta(l_0)>0$ such that
\beqlb\label{new con2}
br+g(y+r)-g(y)\ge \eta(l_0) r
,\quad y\ge0,\quad r\ge l_0.
\eeqlb
Recall that $C_f>0$ is the Lipschitz constant of $f$. Then
\beqlb\label{bound of f}
f(x)\le f(0)+C_fx,\quad x \ge 0.
\eeqlb
\bglemma\label{le:bound of moment}
Suppose that Condition \ref{con2} holds, and let $l_0,\eta(l_0)>0$ be fixed constants that satisfy \eqref{new con2}. Then for all $A\in(0,\eta(l_0)/(2C_f)]$ and $a \in (0, A]$,
$$
\sup_{t}\E [X_t]<\infty,\quad \sup_{t}\E [X_t^2]<\infty,
$$
and both bounds can be chosen uniformly for $a\in(0,A]$.
\edlemma

\proof
 Fix $A\in(0,\frac{\eta(l_0)}{2C_f}]$ and let $a\in(0,A]$. Taking expectations in \eqref{X}, by \eqref{bound of f}, we have 
\beqnn
\E [X_t]\ar=\ar\E [X_0]+\int_0^t\E [-bX_s-g(X_s)]\d s+a\int_0^t\E [f(X_s)]\d s\cr
\ar\le \ar
\E [X_0]+Af(0)t+\int_0^t\E [(-b+AC_f)X_s-g(X_s)]\d s.
\eeqnn
 Then \eqref{new con2} gives 
\beqnn
\E [X_t]\ar\le\ar \E [X_0]+Af(0)t+\int_0^t\E \left[(-b+AC_f)X_s\I_{\{X_s\le l_0\}}\right]\d s\cr
\ar\ar+\int_0^t\E \left[(AC_f-\eta(l_0))X_s\I_{\{X_s>l_0\}}\right]\d s\cr
\ar\le\ar
\E [X_0]+\left[Af(0)+l_0(|b|+AC_f)\right]t-\frac{\eta(l_0)}{2}\int_0^t\E X_s\I_{\{X_s>l_0\}}\d s.
\eeqnn
Using the elementary inequality
\beqlb\label{element ineq}
x\I_{\{x>c\}}\ge x-c,\quad x\ge0,\quad c>0,
\eeqlb
we have
\beqnn
\E [X_t]\ar\le\ar \E [X_0]+\left[Af(0)+l_0(|b|+AC_f)\right]t-\frac{\eta(l_0)}{2}\int_0^t\E [X_s-l_0]\d s\cr
\ar=\ar
\E X_0+K_At-\frac{\eta(l_0)}{2}\int_0^t\E [X_s]\d s,
\eeqnn
where $K_A:=Af(0)+l_0(|b|+AC_f+\eta(l_0)/2)$. An application of Gronwall's inequality yields
$$
\E [X_t]\le 2K_A\eta(l_0)^{-1}+\left[\E [X_0]-2K_A\eta(l_0)^{-1}\right]\e^{-\eta(l_0) t/2}.
$$
It follows  that
\beqlb\label{uniform first moment}
\sup_{a\in(0,A]}\sup_{t\geq0}\E [X_t]<\infty.
\eeqlb
On the other hand, by It\^o's formula,
\beqnn
\E [X_t^2]\ar=\ar \E [X_0^2]+2\int_0^t\E [(-bX_s-g(X_s))X_s]\d s
+2a\int_0^t\E [f(X_s)]\E [X_s]\d s\cr
\ar\ar +2\sigma\int_0^t\E [X_s]\d s.
\eeqnn
By \eqref{uniform first moment} and \eqref{bound of f}, there exists a constant $K'_A>0$ such that
\beqnn
\E [X_t^2]\ar\le\ar \E [X_0^2]+K'_At
+2|b|\int_0^t\E [X_s^2\I_{\{X_s\le l_0\}}]\d s-2\eta(l_0)\int_0^t\E [X_s^2\I_{\{X_s>l_0\}}]\d s\cr
\ar\le\ar
\E [X_0^2]+(K'_A+2|b|l_0^2)t-2\eta(l_0)\int_0^t\E [X_s^2\I_{\{X_s>l_0\}}]\d s.
\eeqnn
Using \eqref{element ineq} and a similar argument, we finish the proof.\qed

Let $l_0>0$ be the fixed constant satisfying \eqref{new con2}. We define
\beqlb\label{function psi}
\psi(r):=
\begin{cases}
c_0r+1-\e^{-c_1r},& 0\le r\le l_0,\\
\psi(l_0)+\frac{\psi'(l_0)}{2}\int_0^{r-l_0}
\left[1+\exp\left(2\frac{\psi''(l_0)}{\psi'(l_0)}s\right)\right]\,\d s,& r>l_0,
\end{cases}
\eeqlb
where
$c_0=c_1\e^{-c_1l_0}$ and $c_1=\frac{2|b|+1}{\sigma}.$
\bglemma\label{le:2.1}
The function $\psi$ defined in \eqref{function psi} satisfies the following.

\item[(i)] $\psi\in C^2(\R_+)$ such that $\psi'\ge0$ and $\psi''\le0$ on $\R_+$;

\item[(ii)] for all $r\ge0$, $
c_0r\le \psi(r)\le (c_0+c_1)r
$ and $\psi'(r)\le c_0+c_1$.
\edlemma

 \proof It is easy to see that
\beqnn
\psi'(r)=
\begin{cases}
c_0+c_1\e^{-c_1r},& 0\le r\le l_0,\\
\frac{\psi'(l_0)}{2}\left[1+\exp\left(2\frac{\psi''(l_0)}{\psi'(l_0)}(r-l_0)\right)\right],& r>l_0,
\end{cases}
\eeqnn
and
\beqnn
\psi''(r)=
\begin{cases}
-c_1^2\e^{-c_1r},& 0\le r\le l_0,\\
\psi''(l_0)\exp\left(2\frac{\psi''(l_0)}{\psi'(l_0)}(r-l_0)\right),& r>l_0.
\end{cases}
\eeqnn
Thus
$\psi\in C^2(\R_+)$, $\psi'>0$, and $\psi''\le0$. This proves assertion (i). 

For $0\leq r\leq l_0$, we have $\psi(r)=c_0r+1-\e^{-c_1r}\geq c_0r$. It together with the fact that $\psi'>0$ implies that $\psi(r)\ge c_0r$ for all $r\ge0$. Moreover, due to the fact that $\psi''\le0$, we have $\sup_{r\ge0}\psi'(r)=\psi'(0)=c_0+c_1$. Hence the assertion (ii) holds.  \qed

We first give an estimate for the operator $\tilde{\mathcal{L}}^i$ given by \eqref{coupling 2.15}. For simplicity, given any function $F\in C^2(\R_+^2)$, we write
\beqlb\label{coupling generator}
\tilde{\mathcal{L}} F(x,y)\ar=\ar
[-bx-g(x)]F_{x}'(x,y)+[-by-g(y)]F_{y}'(x,y)\cr
\ar\ar+\sigma xF''_{xx}(x,y)+\sigma yF_{yy}''(x,y)+2\sigma(x\wedge y)F''_{xy}(x,y).
\eeqlb
Take
$$
F(x,y):=\psi(|x-y|),\quad x,y\in\R_+.
$$

We have the following result.

\bgproposition\label{prop:2.1}
Suppose that Condition \ref{con2} holds. Let $l_0,\eta(l_0)>0$ be fixed constants satisfying \eqref{new con2}. Then
\[
\tilde{\mathcal{L}}F(x,y)\le -\lambda_0 F(x,y),\qquad x\ne y,
\]
where \(\lambda_0=c_0\min\{1,\eta(l_0)\}/(c_0+c_1)\).
\edproposition

\proof It suffices to consider $x>y\ge0$. The case of $0\le x<y$ is symmetric. Note that
\beqlb\label{ineq:2.18}
\tilde{\mathcal{L}} F(x,y)\ar=\ar
[-bx-g(x)]F_{x}'(x,y)+[-by-g(y)]F_{y}'(x,y)\cr
\ar\ar+\sigma xF''_{xx}(x,y)+\sigma yF_{yy}''(x,y)+2\sigma yF''_{xy}(x,y)\cr
\ar=\ar-b(x-y)\psi'(x-y)-[g(x)-g(y)]\psi'(x-y)
+\sigma(x-y)\psi''(x-y).
\eeqlb
It follows that for any $0<x-y\le l_0$,
\beqnn
\tilde{\mathcal{L}} F(x,y)\ar\le\ar|b|(x-y)(c_0+c_1\e^{-c_1(x-y)})
-c_1^2\sigma(x-y)\e^{-c_1(x-y)}\cr
\ar\le\ar
2|b|c_1(x-y)\e^{-c_1(x-y)}-c_1^2\sigma(x-y)\e^{-c_1(x-y)}\cr
\ar\le\ar
-c_1\e^{-c_1(x-y)}(x-y),
\eeqnn
where the last inequality is a consequence of the definition of $c_1$. Then we have
\beqlb\label{ineq:2.19}
\tilde{\mathcal L}F(x,y)\le-c_0(x-y),\quad 0< x-y\le l_0.
\eeqlb
On the other hand, it follows from \eqref{new con2}, \eqref{ineq:2.18} and $\psi''\le0$ that
\beqlb\label{ineq:2.20}
\tilde{\mathcal L}F(x,y)
\le -\eta(l_0)(x-y)\psi'(x-y)\le -\eta(l_0) c_0(x-y)
\eeqlb
for $x - y > l_0$, where the last inequality follows from the fact that $\psi'(r)=c_0(1+\e^{-c_1(r-l_0)})\ge c_0$ for $r>l_0.$
Combining \eqref{ineq:2.19} and \eqref{ineq:2.20}, we have
\[
\tilde{\mathcal L}F(x,y)\le -c_0\min\{1,\eta(l_0)\}|x-y|,\quad x\neq y.
\]
Finally, the assertion follows from Lemma \ref{le:2.1} (ii). \qed

For any $x=(x^1,...,x^N), y=(y^1,...,y^N)$, we write
$$
F_N(x,y):=\frac{1}{N}\sum_{j=1}^N F(x^j,y^j).
$$
Note that $F$ is not differentiable on the diagonal $\Delta=\{(x,y): x=y\}$. Set $D_t^i=\bar X_t^i-X_t^{i,N}$. The quadratic variation of its continuous martingale part satisfies $\d\langle D^{i,c}\rangle_t=2\sigma|D_t^i|\,\d t.$
Consequently, the occupation-density formula gives
\[
L_t^0(D^i)=\lim_{\varepsilon\downarrow0}\frac1{2\varepsilon}
\int_0^t\I_{\{|D_s^i|\le\varepsilon\}}2\sigma|D_s^i|\,\d s=0
\]
by dominated convergence. Thus, the local-time term in the It\^o--Tanaka formula vanishes. For each $1\le i\le N$, applying the It\^{o}-Tanaka formula (see, e.g., Situ \cite[Theorem 142]{Situ}) to $F(\bar{X}^i_t, X^{i,N}_t)$, we obtain
\beqnn
\d F(\bar{X}^i_t, X^{i,N}_t) \ar=\ar \tilde{\mathcal L}F(\bar{X}^i_t, X^{i,N}_t)\I_{\{\bar{X}^i_t\neq X^{i,N}_t\}}\,\d t+ a\E f(X_t)F'_1(\bar{X}^i_t, X^{i,N}_t)\I_{\{\bar{X}^i_t\neq X^{i,N}_t\}}\,\d t \cr
\ar\ar + a\bar{f}^N_t F'_2(\bar{X}^i_t, X^{i,N}_t)\I_{\{\bar{X}^i_t\neq X^{i,N}_t\}}\,\d t +\d M^i_t,
\eeqnn
where $M^i_t$ is a local martingale and $\tilde{\mathcal L}$ is given by \eqref{coupling generator}.
By a standard stopping time argument and Proposition \ref{prop:2.1}, we obtain from summing over $i$ that
\beqlb\label{ineq:2.21}
\ar\ar\E [F_N(\bar{X}_t, X^{N}_t)]-
\E [F_N(\bar{X}_0, X^{N}_0)]\cr
\ar\ar\qquad \le -\lambda_0 \int_0^t\E [F_N(\bar{X}_s, X^{N}_s)]\,\d s\cr
\ar\ar\qquad\qquad +\frac{a}{N}\sum_{i=1}^N
\int_0^t\E\left[[\E [f(X_s)]-\bar f_s^N]\sgn(D_s^i)
\psi'(|D_s^i|)\I_{\{D_s^i\ne0\}}\right]\,\d s\cr
\ar\ar\qquad\le -\lambda_0\int_0^t\E [F_N(\bar{X}_s,X^N_s)]\,\d s+
a(c_0+c_1)\int_0^t\E\left[\left|\E f(X_s)-\bar{f}^N_s\right|\right]\,\d s,
\eeqlb
 where in the last inequality, we have used Lemma \ref{le:2.1} (ii).
For any $A>0$, define
\[
M_A:=\sup_{0<a\leq A}\sup_{t\geq0}(\E[f(X_t)^2])^{1/2}.
\]

\bglemma\label{le:2.3}
 Suppose that Condition \ref{con2} holds. Let $l_0,\eta(l_0)>0$ be fixed constants satisfying \eqref{new con2}. Let $A\in(0,\eta(l_0)/(2C_f)]$. Then $M_A<\infty$ and for every $a\in(0,A]$,
$$
\E\left[\left|\E[f(X_t)]-\bar{f}^N_t\right|\right]
\le \frac{M_A}{\sqrt{N}}+
\frac{C_f}{c_0}\E [F_N(\bar{X}_t,X^N_t)].
$$
\edlemma

\proof The finiteness of $M_A$ follows from
Lemma~\ref{le:bound of moment} and \eqref{bound of f}. Note that
\beqnn
\E [f(X_t)]-\bar{f}^N_t
=\left[\E [f(X_t)]-\frac{1}{N}\sum_{j=1}^N f(\bar{X}^j_t)\right]+
\frac{1}{N}\sum_{j=1}^N
\left[f(\bar{X}^j_t)-f(X^{j,N}_t)\right]:= I_1+I_2.
\eeqnn
Since $\{(\bar{X}^j_t)_{t\ge0}\}_{j=1}^N$ are i.i.d. and $Law(\bar{X}^j_t)=Law(X_t)$ for each $t$, we have $\E[I_1]=0$. It follows from the Cauchy-Schwarz inequality that
$$
\E[|I_1|]\le (\E[|I_1|^2])^{1/2}= \frac{\sqrt{\var(f(\bar{X}^1_t))}}{\sqrt{N}}\le \frac{M_A}{\sqrt{N}}.
$$
On the other hand, by the Lipschitz continuity of $f$ again,
$$
\E[|I_2|]\le \frac{C_f}{N}\sum_{j=1}^N\E\left[|\bar{X}^j_t-X^{j,N}_t|\right].
$$
Thanks to Lemma \ref{le:2.1} (ii), we complete the proof. \qed

 {\it Proof of Theorem \ref{main result}.}
Fix $0<A<a_*$.  By the definition of the supremum in
\eqref{eq:a-star}, there are $l_0,\eta(l_0)>0$ such that
\[
 A<\frac{1}{C_f}\frac{\min\{1,\eta(l_0)\}}{(1+\e^{c_1l_0})^2}<\frac{\eta(l_0)}{4C_f}<\frac{\eta(l_0)}{2C_f}.
\]
Hence Lemma~\ref{le:bound of moment}
and Lemma~\ref{le:2.3} apply uniformly for $a\in(0,A]$.

Write
\[
\lambda:=\lambda_0-\frac{A(c_0+c_1)C_f}{c_0}>0,
\]
where $\lambda_0$ is given in Proposition \ref{prop:2.1}. It follows from \eqref{ineq:2.21} and Lemma~\ref{le:2.3} with this $A$ that,
for every $a\in(0,A]$,
\beqnn
\ar\ar\E [F_N(\bar{X}_t,X^N_t)]-\E [F_N(\bar{X}_0,X^N_0)]\cr
\ar\ar\qquad\le-\left(\lambda_0-\frac{a(c_0+c_1)C_f}{c_0}\right)\int_0^t\E[F_N(\bar{X}_s,X^N_s)]\,\d s+a(c_0+c_1)\int_0^t\frac{M_A}{\sqrt{N}}\,\d s\cr
\ar\ar\qquad \le -\lambda\int_0^t\E [F_N(\bar{X}_s,X^N_s)]\,\d s+A(c_0+c_1)\int_0^t\frac{M_A}{\sqrt{N}}\,\d s.
\eeqnn
By Gronwall's inequality,
\beqlb\label{ineq:2.23}
\E [F_N(\bar{X}_t,X^N_t)]\le \e^{-\lambda t}\E [F_N(\bar{X}_0,X^N_0)]+\frac{A(c_0+c_1)M_A}{\lambda}N^{-1/2}.
\eeqlb
Using Lemma \ref{le:2.1} (ii) gives
$$
\frac{c_0}{N}\sum_{j=1}^N \E\left[|\bar{X}^j_t - X^{j,N}_t|\right]
\le \E[F_N(\bar{X}_t,X^N_t)],
$$
which follows that
\[
W_{d^N}(\mu^{\bigotimes N}_t,\mu^N_t)
\le \frac1{c_0}\E[F_N(\bar{X}_t,
X^N_t)]
\]
since the constructed joint law is a coupling of $\mu_t^{\otimes N}$ and $\mu_t^N$. By the optimal initial coupling \eqref{eq:optimal-initial-coupling} and Lemma \ref{le:2.1} (ii), we get
\[
\E[F_N(\bar{X}_0,X^N_0)]\le (c_0+c_1) W_{d^N}(\mu_0^{\otimes N},(\mu_0^N)^{\otimes N}).
\]
Combining those estimates with \eqref{ineq:2.23} gives the assertion with
$C_1=(c_0+c_1)/c_0$ and $C_2=A(c_0+c_1)M_A/(\lambda c_0)$. This completes the proof of Theorem \ref{main result}.
\qed

{\it Proof of Theorem \ref{propagation 2}.}
We use the coupling constructed in the proof of Theorem~\ref{main result}. Fix $T<\infty$.
By Proposition~\ref{exist and unique} and the Lipschitz continuity of \(f\),
we have \(M_{T,a}<\infty\). The occupation-density argument preceding
\eqref{ineq:2.21} and the monotonicity of \(g\) yield
\beqnn 
\frac1N\sum_{i=1}^N
\E|\bar X_t^i-X_t^{i,N}|
\ar\leq\ar 
\frac1N\sum_{i=1}^N
\E|\bar X_0^i-X_0^{i,N}|   
+b^-\int_0^t\frac1N\sum_{i=1}^N
\E|\bar X_s^i-X_s^{i,N}|\,\d s \cr
\ar\ar 
+a\int_0^t
\E\left|\E[f(X_s)]-\bar f_s^N\right|\,\d s.
\eeqnn
As in the proof of Lemma~\ref{le:2.3},
\[
\E\left|\E[f(X_s)]-\bar f_s^N\right|
\leq M_{T,a}N^{-1/2}
+\frac{C_f}{N}\sum_{i=1}^N
\E|\bar X_s^i-X_s^{i,N}|,
\qquad s\leq T.
\]
Therefore, an application of Gronwall's inequality gives
\[
\frac1N\sum_{i=1}^N\E|\bar X_t^i-X_t^{i,N}|
\leq
\e^{\kappa_at}
\frac1N\sum_{i=1}^N\E|\bar X_0^i-X_0^{i,N}|
+aM_{T,a}\Phi_{\kappa_a}(t)N^{-1/2},
\]
where \(\kappa_a=b^-+aC_f\). The optimality of the initial coupling and
the definition of \(W_{d^N}\) imply the asserted estimate.
\qed

The scaling limit in
Section~\ref{sec:discrete-approximation} also needs convergence of paths on a finite time
interval.  The following proposition is a consequence of the coupling estimate above and
is used for this purpose.

\begin{pro1}\label{prop:pathwise-chaos}
Suppose that Condition~\ref{con2} holds  and $\mu_0 = \mu_0^N$. Fix
$0<A<a_*$ and let $a\in(0,A]$.  Then, for every $T<\infty$, there is a
constant $C_T<\infty$ such that
\[
\E\left[\sup_{0\leq s\leq T}|\bar X_s^1-X_s^{1,N}|\right]
\leq C_TN^{-1/4}.
\]
Consequently,
\[
X^{1,N}\Rightarrow X\qquad\text{in }\D([0,T],\R_+),\qquad N\to\infty.
\]
\end{pro1}

\proof
Choose the identical initial coupling, so $D_0^i=0$.  Estimate
\eqref{ineq:2.23}, exchangeability and Lemma~\ref{le:2.1} imply, with constants
independent of $N$,
\begin{equation}\label{eq:path-coupling-one-time}
\sup_{t\geq0}\E \left[F_N(\bar X_t,X_t^N)\right]\leq C N^{-1/2},
\qquad
\sup_{t\geq0}\E\left[|D_t^1|\right]\leq C N^{-1/2}.
\end{equation}
Lemma~\ref{le:2.3} and the first bound also give
\begin{equation}\label{eq:path-empirical-error}
\sup_{t\geq0}\E\left[|\E \left[f(X_t)\right]-\bar f_t^N|\right]\leq C N^{-1/2}.
\end{equation}

Apply the Tanaka formula to $|D_t^1|$.  Its local time at zero vanishes by
\eqref{coupling 2.15} and the occupation-density calculation preceding
\eqref{ineq:2.21}; moreover, the contribution of $g$ is non-positive. Hence,
\[
|D_t^1|
\leq M_t+b^-\int_0^t|D_s^1|\,\d s
+a\int_0^t|\E f(X_s)-\bar f_s^N|\,\d s,
\]
where $M$ is a continuous local martingale satisfying $\d\langle M\rangle_t\leq2\sigma|D_t^1|\,\d t.$
After localization, by the Burkholder--Davis--Gundy inequality,
\eqref{eq:path-coupling-one-time} and \eqref{eq:path-empirical-error}, we have
\beqnn 
\E\left[\sup_{0\leq t\leq T}|D_t^1|\right]
\ar\leq\ar  C\left(\int_0^T\E\left[|D_s^1|\right]\,\d s\right)^{1/2}
+b^-\int_0^T\E\left[|D_s^1|\right]\,\d s\cr
\ar\ar +a\int_0^T\E\left[|\E \left[f(X_s)\right]-\bar f_s^N|\right]\,\d s
\leq C_TN^{-1/4}.
\eeqnn
Fatou's lemma removes the localization. Thus
$X^{1,N}-\bar X^1\to0$ uniformly on $[0,T]$ in probability.  Since $\bar X^1$
has the same path law as $X$, the asserted convergence follows. \qed

\section{Proof of Theorem~\ref{thm:discrete-approximation}}\label{sec:discrete-approximation}

This section proves Theorem~\ref{thm:discrete-approximation}. We use convergence of
generators and martingale problems; see Ethier and Kurtz \cite[Chapter~3, Sections~5 and~8, and Chapter~4,
Section~8]{EK86}. For $N\geq1$, set $\mathcal D_N:=\big\{F|_{\R_+^N}:F\in C_c^3(\R^N)\big\}.$
Thus, each $F\in\mathcal D_N$ has a $C^3$ extension to an open neighborhood of $\R_+^N$;
Taylor expansions at the boundary are understood through such an extension.

Define
\[
f^N(x)=\frac1N\sum_{j=1}^Nf(x^j),\qquad x\in\R_+^N,
\]
and the operator
\begin{equation}\label{eq:continuous-particle-generator-scaling}
\mathcal A_NF(x)=\sum_{i=1}^N\left\{\sigma x^i\partial_{ii}F(x)
+[-bx^i-g(x^i)+af^N(x)]\partial_iF(x)\right\},\qquad F\in\mathcal D_N.
\end{equation}
This is the generator of the continuous particle system \eqref{Xi}. The following result follows directly from Proposition \ref{prop:particle-wellposedness}; its proof is omitted here.

{
\begin{lem1}\label{lem:continuous-martingale-problem}
For every fixed $N$ and every deterministic initial point in $\R_+^N$, the martingale
problem for $\mathcal A_N$ on $\mathcal D_N$ is well posed.
\end{lem1}

}

\begin{lem1}\label{lem:discrete-nonexplosion}
For every fixed $n,N\geq1$ and every deterministic initial state
$z\in\mathbb Z_+^N$, the Markov chain $Z^{n,N}$ with rates
\eqref{eq:discrete-birth-rate}--\eqref{eq:discrete-death-rate} is non-explosive.
\end{lem1}

\proof
Let $|z|_1:=\sum_{i=1}^N z_i$ and $S_t^{n,N}:=\sum_{i=1}^N Z_t^{i,n,N}.$
At state $z$, the total rate of upward jumps is
\beqnn
\Lambda^{+,n,N}(z)
:=\sum_{i=1}^N\lambda_i^{+,n,N}(z) =(n\sigma+b^-)|z|_1
  +na\sum_{j=1}^N f(z_j/n).
\eeqnn
By \eqref{bound of f}, we have $\Lambda^{+,n,N}(z)
\le c_2+c_3|z|_1,$
where $c_2:=naNf(0)$ and $c_3:=n\sigma+b^-+aC_f.$

Let \(s_0:=S_0^{n,N}\), and let \(\Pi(\mathrm ds,\mathrm du)\) be a
Poisson random measure on \(\mathbb R_+^2\) with intensity
\(\mathrm ds\,\mathrm du\). Define
\[
Y_t
:=
s_0+
\int_{(0,t]}\int_0^\infty
\mathbf 1_{\{u\le c_2+c_3Y_{s-}\}}
\,\Pi(\mathrm ds,\mathrm du).
\]
Thus, \(Y\) is the linear pure-birth process with immigration having
an initial state \(Y_0=s_0\) and transition rates
\[
k\longrightarrow k+1
\quad\text{at rate}\quad
c_2+c_3k,\qquad k\in\mathbb Z_+.
\]
Here, \(c_2\) is the immigration rate and \(c_3\) is the
 birth rate per-particle. The process \(Y\) is nonexplosive. Indeed, if
\(c_2+c_3s_0>0\), then $
\sum_{k=s_0}^{\infty}(c_2+c_3k)^{-1}=\infty.$
Hence, by the standard nonexplosion criterion for pure-birth processes
(see, e.g., Norris~\cite[Theorem~2.5.2]{Norris97}), the lifetime of
\(Y\) is almost surely infinite. If \(c_2+c_3s_0=0\), then, since
\(c_2\geq0\) and \(c_3>0\), we necessarily have \(c_2=s_0=0\).
Consequently, \(Y_t=0\) for every \(t\geq0\). Thus \(Y\) is
nonexplosive in either case.

On $E:=\{(z,y)\in\mathbb Z_+^N\times\mathbb Z_+: |z|_1\le y\},$
we couple \(Z^{n,N}\) and \(Y\) as follows:
\[
(z,y)\longrightarrow(z+e_i,y+1)
\quad\text{at rate}\quad
\lambda_i^{+,n,N}(z),
\]
\[
(z,y)\longrightarrow(z,y+1)
\quad\text{at rate}\quad
c_2+c_3y-\Lambda^{+,n,N}(z),
\]
and
\[
(z,y)\longrightarrow(z-e_i,y)
\quad\text{at rate}\quad
\lambda_i^{-,n,N}(z).
\]
Notice that $\Lambda^{+,n,N}(z)\le c_2+c_3|z|_1\le c_2+c_3y$ on $E$.
Moreover, $E$ is invariant under these transitions. Consequently, $S_t^{n,N}\le Y_t$
up to the lifetime of the coupled process, and every upward jump of
$Z^{n,N}$ is accompanied by an upward jump of $Y$. The result follows from the non-explosiveness of $Y$.
\qed

By the definition of $X^{n,N}$ and \eqref{eq:discrete-birth-rate}-\eqref{eq:discrete-death-rate}, one sees that the generator $\mathcal A_{n,N}$ of $X^{n,N}$ acts on $F\in\mathcal D_N$ as
\beqlb\label{eq:discrete-generator}
\mathcal A_{n,N}F(x)
\ar=\ar \sum_{i=1}^N\big(n^2\sigma x^i+nb^-x^i+naf^N(x)\big)
   [F(x+n^{-1}e_i)-F(x)]\cr
\ar\ar +\sum_{i=1}^N\I_{\{x^i\geq n^{-1}\}}
  \big(n^2\sigma x^i+nb^+x^i+ng(x^i)\big)
  [F(x-n^{-1}e_i)-F(x)].
\eeqlb

\begin{lem1}\label{generator_convergence}
 Let $\mathcal A_{n,N}$ and $\mathcal{A}_N$ be the generators of $X^{n,N}$ and $X^N$ given by \eqref{eq:discrete-generator} and \eqref{eq:continuous-particle-generator-scaling}, respectively. For any $F\in\mathcal D_N$, we have
 \beqnn
\lim_{n\to\infty}\sup_{x\in n^{-1}\mathbb Z_+^N}
|\mathcal A_{n,N}F(x)-\mathcal A_NF(x)|=0.
\eeqnn
\end{lem1}
\proof
Let $e_i$ be the $i$th coordinate vector of $\R^N$. Uniformly on compact subsets of $\R_+^N$, as $n \rightarrow \infty,$ we have
\beqnn
F(x \pm n^{-1}e_i)-F(x)
=\pm n^{-1}\partial_iF(x)+\tfrac12n^{-2}\partial_{ii}F(x)+O(n^{-3}).
\eeqnn
By the above, \eqref{eq:continuous-particle-generator-scaling} and \eqref{eq:discrete-generator}, for any $F\in\mathcal D_N$, we have
\beqnn
\mathcal A_{n,N}F(x)
  \ar=\ar \sum_{i = 1}^N (-bx^i - g(x^i) + af^N(x))\partial_i F(x) + \sum_{i = 1}^N (\sigma x^i + O(n^{-1}))\partial_{ii}F(x) + O(n^{-1})\cr
  \ar=\ar \mathcal{A}_N F(x) + O(n^{-1})\sum_{i = 1}^N \partial_{ii}F(x) + O(n^{-1})
\eeqnn
for $x\in n^{-1}\mathbb Z_+^N$.
Recall that $F$ and its derivatives are compactly
supported. Moreover, $f$ and $g$ are bounded in this compact set. The result follows.
\qed

\begin{lem1}\label{tight}
The sequence $\{X^{n,N}: n \ge 1\}$ is tight in $D([0, T], \mathbb{R}_+^N)$.
\end{lem1}

\proof
For each coordinate of  $X^{n,N}$, we have
\begin{equation}\label{eq:discrete-semimartingale}
X_t^{i,n,N}=X_0^{i,n,N}+M_t^{i,n,N}
+\int_0^t[-bX_s^{i,n,N}-g(X_s^{i,n,N})+af^N(X_s^{n,N})]\,\d s,
\end{equation}
where $M^{i,n,N}$ is a local martingale with predictable quadratic variation
\beqnn
\langle M^{i,n,N}\rangle_t
=\int_0^t\left[2\sigma X_s^{i,n,N}
+\frac{|b|X_s^{i,n,N}+g(X_s^{i,n,N})+af^N(X_s^{n,N})}{n}\right]\d s.
\eeqnn
Let $V(x)=1+\sum_{i=1}^Nx^i$.  By \eqref{eq:discrete-generator}, \eqref{bound of f} and $g\geq0$, one obtains
\beqlb\label{eq:discrete-lyapunov}
\mathcal A_{n,N}V(x)
=-b\sum_{i=1}^Nx^i-\sum_{i=1}^Ng(x^i)+a\sum_{i=1}^Nf(x^i)
\leq C_NV(x)
\eeqlb
with $C_N$ independent of $n$.  Moreover, we have
$\sup_n\E [V(X_0^{n,N})]<\infty$ since $X_0^{i,n,N}\leq\xi_i$.
For $m>1$, we define $\tau_m^{n,N}=\inf\{t\geq0:V(X_t^{n,N})\geq m\}.$
By \eqref{eq:discrete-lyapunov} and It\^{o}'s formula, it follows that
\beqnn
\E [V(X_{t\wedge\tau_m^{n,N}}^{n,N})] \ar=\ar \E[V(X_0^{n,N})] + \E\left[\int_0^{t\wedge \tau_m^{n,N}} \mathcal A_{n,N}V(X_s^{n,N})\,\d s\right]\\
\ar\leq\ar \E [V(X_0^{n,N})]+C_N\int_0^t
\E [V(X_{s\wedge\tau_m^{n,N}}^{n,N})]\,\d s.
\eeqnn
By Gronwall's inequality, there exists a constant $C_{T,N}$ independent of $m$ and $n$ such that $\sup_{n\geq1}\E [V(X_{T\wedge\tau_m^{n,N}}^{n,N})] \leq C_{T,N}.$
Consequently,
\beqlb\label{tau}
\sup_{n\geq1}\pP(\tau_m^{n,N}\leq T)\leq\frac{C_{T,N}}m\longrightarrow0,
\qquad m\to\infty.
\eeqlb
Fix $m\geq1$. By \eqref{eq:discrete-semimartingale} define
\beqnn
B^{i,n,N,m}_t := \int_0^{t\wedge \tau_m^{n,N}}[-bX_s^{i,n,N}-g(X_s^{i,n,N})+af^N(X_s^{n,N})]\,\d s
\eeqnn
and $M_t^{i,n,N,m} = M_{t\wedge\tau_m^{n,N}}^{i,n,N}$. Let $B^{n,N,m}$ and $M^{n,N,m}$ be the corresponding vectors. Here $|z| := (\sum_{i = 1}^N (z^i)^2)^{1/2}$ denotes the Euclidean norm of $z \in \mathbb{R}^N$. On \(\{\tau_m^{n,N}>0\}\), the process remains in the region
\(\{V\leq m\}\) before the stopping time. By the overshoot estimate, $V(X_{\tau_m^{n,N}}^{n,N})
\leq m+\frac1n\leq m+1.$
On $\{\tau_m^{n,N}=0\}$, the stopped process is
constant, and hence both $B^{n,N,m}$ and $M^{n,N,m}$ have zero
increments. One can easily check that $M^{n,N,m}$ is square-integrable.

Recall that
\(X_0^{i,n,N}\leq\xi_i\) for \(1\leq i\leq N\), and hence $
\sup_{n\geq1}
\mathbb P\bigl(|X_0^{n,N}|>R\bigr)
\rightarrow 0$ as $R \rightarrow \infty.$
Moreover, by Doob's inequality, we have
\beqnn
\mathbb E\left[
 \sup_{t\leq T}|M_t^{n,N,m}|^2
\right]
\leq
4\mathbb E[|M_T^{n,N,m}|^2]
=
4\sum_{i = 1}^N \mathbb{E} [\langle M^{i, n,N,m}\rangle_T]
\leq4C_{m,N}T.
\eeqnn
Therefore,
\begin{align*}
\mathbb P\left(
 \sup_{t\leq T}|X_{t\wedge\tau_m^{n,N}}^{n,N}|>R
\right)
&\leq
\mathbb P\left(|X_0^{n,N}|>\frac R3\right)
 +\mathbb P\left(
  \operatorname{Var}_{[0,T]}(B^{n,N,m})>\frac R3
 \right)\\
&\quad
 +\mathbb P\left(
  \sup_{t\leq T}|M_t^{n,N,m}|>\frac R3
 \right)\\
&\leq
\mathbb P\left(|X_0^{n,N}|>\frac R3\right)
+\frac{3C_{m,N}T}{R}
+\frac{36C_{m,N}T}{R^2},
\end{align*}
where $\operatorname{Var}_{[0,T]}(B^{n,N,m})$ is the total variation on $[0,T]$ of $B^{n,N,m}$. The right-hand side converges to zero as \(R\to\infty\), uniformly in
\(n\). Thus the stopped processes satisfy compact containment in
\(\mathbb R_+^N\).

We next verify the Aldous modulus condition. Let
\(\rho_n\leq T\) be stopping times, let \(\theta_n\downarrow0\), and set $\sigma_n:=(\rho_n+\theta_n)\wedge T.$
It follows that $\mathbb E[|B_{\sigma_n}^{n,N,m}-B_{\rho_n}^{n,N,m}|]\leq C_{m,N}\theta_n$ and
\beqnn
\mathbb E\left[\left|
 M_{\sigma_n}^{n,N,m}-M_{\rho_n}^{n,N,m}
\right|^2\right]
 =
\sum_{i = 1}^N\mathbb E \left[
 \langle M^{i,n,N,m}\rangle_{\sigma_n}
 -
 \langle M^{i,n,N,m}\rangle_{\rho_n}
\right]
 \leq C_{m,N}\theta_n.
\eeqnn
Markov's inequality yields, for every \(\varepsilon>0\),
\beqnn
\mathbb P\left(
 \left|X_{\sigma_n\wedge\tau_m^{n,N}}^{n,N}
       -X_{\rho_n\wedge\tau_m^{n,N}}^{n,N}\right|>\varepsilon
\right)
\ar\leq\ar
\frac{2}{\varepsilon}
\mathbb E\left[\left|
 B_{\sigma_n}^{n,N,m}-B_{\rho_n}^{n,N,m}
\right|\right] +
\frac{4}{\varepsilon^2}
\mathbb E\left[\left|
 M_{\sigma_n}^{n,N,m}-M_{\rho_n}^{n,N,m}
\right|^2\right]\cr
\ar\leq\ar
C_{m,N}\left(
 \frac{2}{\varepsilon}+\frac{4}{\varepsilon^2}
\right)\theta_n
\rightarrow 0
\eeqnn
as $n \rightarrow \infty$, which gives Aldous'
criterion \cite{Aldous78}; see also \cite[Theorem~3.8.6]{EK86}.
Then for every fixed $m$, $\{
X_{\cdot\wedge\tau_m^{n,N}}^{n,N}:n\geq1\}$
is tight in \(\mathbb D([0,T],\mathbb R_+^N)\). By \eqref{tau}, it follows that $\{X^{n,N}:n\geq1\}$
is tight in \(\mathbb D([0,T],\mathbb R_+^N)\).
\qed

{\it Proof of Theorem~\ref{thm:discrete-approximation}}.
We now identify the sub-sequential limits for fixed $N$.
By Lemma \ref{tight} and Prokhorov's theorem, there exists a convergent subsequence; see \cite[Theorems~3.2.2 and~3.5.6]{EK86}.  Let $Y$ be its limit.  For every
$F\in\mathcal D_N$,
\beqnn
F(X_t^{n,N})-F(X_0^{n,N})-
\int_0^t\mathcal A_{n,N}F(X_s^{n,N})\,\d s
\eeqnn
is a martingale.  Recall that $X_0^{n,N}\to(\xi_1,\ldots,\xi_N)$ almost surely. Let $0\leq t_1<\cdots<t_k\leq s<t\leq T$ be continuity times of $Y$ and let
$\Phi$ be a bounded and continuous function. Then by the standard martingale-problem convergence argument (see, e.g.,\cite[Theorem~4.8.10]{EK86}) and Lemma \ref{generator_convergence}, we obtain 
\beqnn
\E\left[\Phi(Y_{t_1},\ldots,Y_{t_k})
\left(F(Y_t)-F(Y_s)-\int_s^t\mathcal A_NF(Y_r)\,\d r\right)\right]=0.
\eeqnn
By right-continuity, one can check that the above holds for all $s<t$.  Hence $Y$ solves the
$\mathcal A_N$-martingale problem on $\mathcal D_N$. By
Lemma~\ref{lem:continuous-martingale-problem}, its law is that of the solution $X^N$ of
\eqref{Xi}. Every subsequential limit is therefore the same, and $X^{n,N}\Rightarrow X^N$ in $\D([0,T],\R_+^N).$
This proves part~(i).

 By Proposition~\ref{prop:pathwise-chaos}, one sees that $X^{1,N}\Rightarrow X$ in
$\D([0,T],\R_+)$. On the other hand, for every bounded continuous
functional $\Phi$ on $\D([0,T],\R_+)$, we have
\beqnn
|\E[\Phi(X^{1,n,N})]-\E[\Phi(X)]|
\leq |\E[\Phi(X^{1,n,N})]-\E[\Phi(X^{1,N})]| +|\E[\Phi(X^{1,N})]-\E[\Phi(X)]|.
\eeqnn
For fixed $N$, the first term tends to zero as $n\to\infty$ by part~(i), and the second
{ tends to zero as $N\to\infty$ by
Proposition~\ref{prop:pathwise-chaos}.}  This proves part~(ii).

For part~(iii), let $d_T$ be the Prokhorov metric in the space of probability measures on
$\D([0,T],\R_+)$.  It metrizes weak convergence because the Skorokhod space is Polish
\cite[Theorems~3.3.1 and~3.5.6]{EK86}.  Choose a strictly increasing sequence $N_k\to\infty$ such that
\[
d_T(\operatorname{Law}(X^{1,N_k}),\operatorname{Law}(X))\leq\frac1{2k}.
\]
For each $k$, part~(i) gives $n_k$ such that, for every $n\geq n_k$,
\[
d_T(\operatorname{Law}(X^{1,n,N_k}),\operatorname{Law}(X^{1,N_k}))\leq\frac1{2k}.
\]
Choose $n_k$ strictly increasing, set $N_T(n)=N_1$ for $n<n_1$, and set
$N_T(n)=N_k$ when $n_k\leq n<n_{k+1}$.  Then $N_T(n)\to\infty$ and, on the latter interval,
\[
d_T(\operatorname{Law}(X^{1,n,N_T(n)}),\operatorname{Law}(X))\leq\frac1k.
\]
This proves part~(iii). \qed

\section{Proofs of Propositions \ref{prop:logistic-threshold}, \ref{prop: bound of ac}, \ref{p0814} and Theorem \ref{thm:logistic-boundary}}\label{example}

In this section, we prove the results for the logistic mean-field branching
diffusion introduced in Subsection~\ref{subsec:logistic-model}. Since the upper
bound in Proposition~\ref{prop: bound of ac} uses the supercritical invariant
law from Proposition~\ref{prop:logistic-threshold}, we first establish the
phase transition and then return to the uniform-in-time threshold. We finally
derive the long-time mean limit needed for
Theorem~\ref{thm:logistic-boundary} and prove the boundary classification by
comparison with frozen-immigration diffusions. Boundary behavior for related
square-root diffusions is discussed by Karlin and Taylor~\cite[Chapter~15]{KT81} and by Peskir
and Roodman~\cite{PR23}.  A more general inhomogeneous framework is
presented in Chen et al.~\cite{CYZ26}.

We begin with a moment estimate used in both the phase-transition and
long-time convergence arguments.

\begin{lem1}\label{lem:logistic-moments}
The solution of \eqref{eq:logistic-mv} satisfies $\sup_{t \ge 0}\mathbb{E}[X_t] + \sup_{t\ge0}\E [X_t^2]<\infty.$
Consequently, $m_t := \mathbb{E}[X_t]$ is uniformly continuous on $[0,\infty).$
\end{lem1}

\proof 
For the logistic coefficients, $C_f=1$ and
$\eta(l)=\beta+b_2l\to\infty$ as $l\to\infty$. Given $a>0$, choose
$l_0>0$ such that $\eta(l_0)\geq2a$. By taking $A = a$ in Lemma~\ref{le:bound of moment},
it yields $\sup_{t\geq0}\mathbb E[X_t]
+\sup_{t\geq0}\mathbb E[X_t^2]<\infty.$
Furthermore, taking expectations in \eqref{eq:logistic-mv}, then for
$0\leq s<t$, it gives $m_t-m_s
=\int_s^t\left((a-\beta)m_r
-b_2\mathbb E[X_r^2]\right)\d r.$
The preceding moment bounds therefore imply
$|m_t-m_s|\leq C|t-s|$ for some constant $C<\infty$. Hence $m$ is
globally Lipschitz and, in particular, uniformly continuous.
\qed

{\it Proof of Proposition \ref{prop:logistic-threshold}}.
Define the Lyapunov function
\beqnn
 r(x)=\frac1{\sigma w(x)}\int_x^\infty w(y)\,\d y\quad \text{and} \quad
 V(x)=\int_0^x r(y)\,\d y,
\eeqnn
where $w(\cdot)$ is given by (\ref{scale fun1}).
One can check that  $r$ is strictly decreasing with $ r(x) \le r(0)=1/a_\star^{\rm log}$,  $r(x)\leq C/(1+x)$ and $V(x)\leq C\bigl(1+\log(1+x)\bigr).$
Moreover, there exists $\kappa>0$ such that $r(0)-r(x)\geq\kappa x/(1+x)$ for $x\geq0.$
Applying It\^o's formula to $V(X_t)$, we obtain
\beqnn
 \mathbb E[V(X_t)] \ar=\ar \mathbb E[V(X_0)]
 + \int_0^t m_s\bigl(a\mathbb E[r(X_s)]-1\bigr)\,\d s\cr
 \ar\le\ar \mathbb E[V(X_0)] + \left(\frac{a}{a_\star^{\rm log}} - 1\right)\int_0^t m_s \d s.
\eeqnn
If $a < a_\star^{\rm log}$, we have $\int_0^\infty m_s < \infty$ by the above. It follows from Lemma \ref{lem:logistic-moments} that $m_t\to0$. If $a=a_\star^{\rm log}$, the preceding identity and the bound above give
\beqnn
\mathbb{E}[V(X_t)] \ar=\ar \mathbb{E}[V(X_0)] + a_\star^{\rm log}\int_0^t m_s(\mathbb{E}[r(X_s)] - r(0)) \d s \cr
\ar\le\ar \mathbb{E}[V(X_0)] - a_\star^{\rm log}\kappa\int_0^t m_s \mathbb{E}\left[\frac{X_s}{1 + X_s}\right]\d s.
\eeqnn
By the Cauchy-Schwarz inequality,
\beqnn
m_t^2=(\mathbb E [X_t])^2
=\left(\mathbb E\left[\sqrt{\frac{X_t}{1+X_t}}\sqrt{X_t(1+X_t)}\right]\right)^2
\le \mathbb E\left[\frac{X_t}{1+X_t}\right]\cdot \mathbb E \left[X_t(1+X_t)\right].
\eeqnn
By Lemma \ref{lem:logistic-moments}, it follows that $\mathbb{E}[V(X_t)] \le \mathbb{E}[V(X_0)] - C\int_0^t m_s^3\d s,$
which implies $\int_0^\infty m_t^3 \d t<\infty$, and the uniform continuity in Lemma \ref{lem:logistic-moments} again yields $m_t\to0$. Since
$X_t\geq0$, this proves~(i).

We now consider the case $a>a_\star^{\rm log}$.  For $m>0$, consider
\beqlb\label{barY}
  \bar{Y}_t
  = \bar{Y}_0 + \int_0^t \bigl(am-\beta \bar{Y}_s-b_2(\bar{Y}_s)^2\bigr)\d s
  + \int_0^t \int_0^{\bar{Y}_s}W(\d s, \d u).
\eeqlb
A classical result (see, e.g., \cite[p. 221]{KT81}) implies that $\bar{Y}$ admits a unique invariant probability measure $\pi_\alpha(\d x)  =  \frac{x^{\alpha-1}w(x)}{I_\alpha} \d x,$
where $\alpha=am/\sigma$, the functions $w$ and $I$ are given by \eqref{scale fun1}. The mean of $\pi_\alpha$ is $I_{\alpha + 1}/I_\alpha$. Then $m = I_{\alpha + 1}/I_\alpha$ if and only if $a = \sigma\alpha I_\alpha/I_{\alpha + 1}$. Let
\beqlb\label{F}
F(x) := \frac{\sigma x I_x}{I_{x+1}}, \qquad x > 0.
\eeqlb
One can easily check that $F$ is strictly increasing from $a_\star^{\rm log}$ to infinity. Then there exists a unique $\alpha_a > 0$ such that
\beqlb\label{eq:alpha-a}
  a = F(\alpha_a) = \frac{\sigma\alpha_a I_{\alpha_a}}
       {I_{\alpha_a+1}}.
\eeqlb
Define $\pi_a(\d x) =  \frac{x^{\alpha_a-1}w(x)}{I_{\alpha_a}}\d x$ and $m_a = \int_0^\infty x\pi_a(\d x) = \frac{I_{\alpha_a + 1}}{I_{\alpha_a}}.$
By \eqref{eq:alpha-a}, $m_a = \sigma\alpha_a/a$. Let $Y$ be a solution of \eqref{barY} with $m = m_a$ and $Law(Y) = \pi_a$. Therefore, $Y$ is a weak solution of \eqref{eq:logistic-mv}. By the weak well-posedness of \eqref{eq:logistic-mv}, the solution of the McKean-Vlasov equation with initial law $\pi_a$ has law $\pi_a$ at every time. Consequently, $\pi_a$ is a nonzero  invariant probability measure of \eqref{eq:logistic-mv}. Conversely, any nonzero invariant law of \eqref{eq:logistic-mv} with finite mean $m$ must be the invariant law of the corresponding frozen diffusion and must satisfy \eqref{eq:alpha-a}. The strict monotonicity of $F$ proves uniqueness and completes (ii).
\qed

{\it Proof of Proposition \ref{prop: bound of ac}}. The logistic coefficients satisfy Condition~\ref{con2}, so
Corollary~\ref{corollary} gives $a_*\leq a_c$. It remains to prove the
upper bound. Let $A>a_\star^{\rm log}$. Suppose that \(A\) is admissible.
Choose an interaction strength $a$ such that $a_\star^{\rm log}<a\leq A.$
Proposition~\ref{prop:logistic-threshold} (ii) gives a nonzero invariant law
$\pi_a$ for \eqref{eq:logistic-mv}. Let $m_a:=\int_0^\infty x\,\pi_a(\d x)>0.$
Initialize the limiting equation with law $\pi_a$ and the particle system with
product law $\pi_a^{\otimes N}$. The limiting product law then
remains $\pi_a^{\otimes N}$ for all time. For the particle system, we set $
S_t^N:=\sum_{i=1}^N X_t^{i,N}$
and define the continuous local martingale
\[
M_t^N:=\sum_{i=1}^N\int_0^t\int_0^{X_s^{i,N}}W^i(\d s,\d u).
\]
The independence of the white noises gives $\langle M^N\rangle_t=2\sigma\int_0^t S_s^N\,\d s.$
By L\'evy's characterization, after extending the probability space if
necessary, there is a Brownian motion $B^N$ such that $M_t^N=\int_0^t\sqrt{2\sigma S_s^N}\,\d B_s^N.$
Consequently, summing the particle equations gives
\begin{equation}\label{eq:prop18-total-mass}
S_t^N=S_0^N+\int_0^t\sqrt{2\sigma S_s^N}\,\d B_s^N
+\int_0^t\left[(a-\beta)S_s^N
-b_2\sum_{i=1}^N(X_s^{i,N})^2\right]\d s.
\end{equation}
Let $Y^N$ be driven by
the same Brownian motion $B^N$ and solve
\begin{equation}\label{eq:prop18-dominating-diffusion}
Y_t^N=S_0^N+\int_0^t\sqrt{2\sigma Y_s^N}\,\d B_s^N
+\int_0^t\left[(a-\beta)Y_s^N
-\frac{b_2}{N}(Y_s^N)^2\right]\d s.
\end{equation}
By the Cauchy--Schwarz inequality,
$\sum_{i=1}^N(X_s^{i,N})^2\geq N^{-1}(S_s^N)^2$. Then the drift of $S^N$ in \eqref{eq:prop18-total-mass} is bounded above by the
logistic drift in \eqref{eq:prop18-dominating-diffusion}. An application of comparison principle (see, e.g., Fu and Li \cite[Theorem 5.4]{FuLi2010}) gives
\begin{equation}\label{eq:prop18-comparison}
S_t^N\leq Y_t^N\quad\text{for all }t\geq0\quad\text{almost surely}.
\end{equation}
Note that the process $Y^N$ is
absorbed at zero in finite time almost surely; see
\cite{Lambert05}. It follows from \eqref{eq:prop18-comparison} that \((S^N_t)_{t\ge0}\) is absorbed at zero in finite time almost surely.

It remains to justify convergence of the first moments. Following a similar argument as in the proof of Lemma \ref{le:2.1}, one has
$\sup_{t\geq0}\mathbb E[(Y_t^N)^2]<\infty$ for every fixed $N$. It implies that $(Y_t^N)_{t\geq0}$ is uniformly integrable and hence
$\mathbb E[Y_t^N]\to0$ as \(t\to\infty\) since $Y_t^N\to0$ as \(t\to\infty\) almost surely. By \eqref{eq:prop18-comparison},
$\mathbb E[S_t^N]\to0$ as \(t\to\infty\). 
Consequently,
\[
W_{d^N}(\mu_t^N,\delta_{\mathbf0})
=\frac1N\mathbb E[S_t^N]\longrightarrow0.
\]
Since \(
W_{d^N}(\pi_a^{\otimes N},\delta_{\mathbf0})=m_a.\)
The triangle inequality therefore gives
\[
\liminf_{t\to\infty}
W_{d^N}(\pi_a^{\otimes N},\mu_t^N)\geq m_a
\qquad\text{for every fixed }N\ge2.
\]
The admissibility of $A$ implies
\[
\lim_{t\to\infty}W_{d^N}(\pi_a^{\otimes N},\mu_t^N)\leq C_2N^{-1/2},
\quad N\geq2.
\]
This is a contradiction since $m_a>0$. Hence $a_c\leq a_\star^{\rm log}$.\qed

{\it Proof of Proposition \ref{p0814}. }
Put $\widehat X_t:=X_{t/a}$. Then $\widehat X$ is the logistic
McKean--Vlasov diffusion in~\cite{Hutz12} with parameters $\gamma=\frac{b_2}{a},  K=\frac{a-\beta}{b_2}, \beta_{\rm H}=\frac{\sigma}{a}.$
These parameters are positive. Indeed, if $\beta>0$, then
$I_1<\sigma/\beta$ and hence $a_\star^{\rm log}>\beta$, while the assertion is
immediate if $\beta\leq0$. Moreover, the left-hand side of condition \cite[(3.10)]{Hutz12} equals
\beqnn
\int_0^\infty
\exp\left(-\frac{\beta z}{a}
-\frac{b_2\sigma z^2}{2a^2}\right)\d z
=\frac{a}{\sigma}I_1
=\frac{a}{a_\star^{\rm log}}>1.
\eeqnn
Thus condition \cite[(3.10)]{Hutz12}  fails. Corollary~3.5 and the argument
in the proof of Corollary~3.6 of \cite{Hutz12} imply that, for every
$x>0$, $\widehat X_t$ converges in distribution to a nontrivial
invariant law. By Lemma~\ref{lem:logistic-moments}, this limiting law has
finite second moment. Then $\pi_a$ is given by Proposition~\ref{prop:logistic-threshold}(ii). Since $X_t=\widehat X_{at}$, it follows that $X_t \Rightarrow \pi_a$ as $t \rightarrow \infty$.

Moreover, by Lemma~\ref{lem:logistic-moments},
\beqnn
\sup_{t\geq0}\mathbb E_x\left[
X_t\mathbf 1_{\{X_t>R\}}\right]
\leq \frac1R\sup_{t\geq0}\mathbb E_x[X_t^2]
\longrightarrow0
\qquad\text{as }R\to\infty.
\eeqnn
Thus $\{X_t:t\geq0\}$ is uniformly integrable, and convergence of
the means follows. The displayed formula for $m_a$ follows from the
definition of $\pi_a$ and \eqref{eq:alpha-a}.
\qed

{\it Proof of Theorem \ref{thm:logistic-boundary}}.
We first classify the boundary at zero for diffusions with a constant
immigration rate. Parts~(i)--(ii) follow by comparing the mean-field solution
with these frozen diffusions, while part~(iii) additionally uses the
supercritical mean limit in Proposition~\ref{p0814}.
For a constant immigration rate $c\geq0$, consider
\beqnn
 Y_t^c= Y_0^c + \int_0^t
 (c-\beta Y_s^c-b_2(Y^c_s)^2)\d s  +\int_0^t\int_0^{Y_s^c}W(\d s, \d u).
\eeqnn
Let $\tau_0^{Y^c} := \inf\{t \ge 0: Y_t^c = 0\}.$
The scale density is $s_c'(x) = \frac{x^{-c/\sigma}}{w(x)}$, where $w(x)$ is given by \eqref{scale fun1}. Notice that $s_c(x) \rightarrow \infty$ as $x \rightarrow \infty$. By Feller's boundary test (see, e.g.,
\cite[Chapter~15, Section~6, especially Table~15.6.2 and
Lemma~15.6.3]{KT81}), zero is accessible precisely when $c < \sigma$. By the two-sided scale formula and $s_c(\infty) = \infty$, for any $x > 0$ we have $\mathbb{P}_x(\tau_0^{Y^c} < \infty) = 0$ when $c \ge \sigma$. Moreover, when $c < \sigma$, we have $\mathbb{P}_x(\tau_0^{Y^c} < \infty) = 1$. By standard finite-time properties of regular one-dimensional diffusions (see, e.g., \cite[Chapter 15, Section 6]{KT81}), it gives, for any $t > 0$ and $x > 0$,
\beqnn
 p_t^c(x) := \mathbb{P}_x(\tau_0^{Y^c} < t) > 0  \qquad (c < \sigma).
\eeqnn

We first prove~(i). Suppose that $am_t\geq\sigma$ for every $t \ge 0$.
Let $Y^\sigma$ be the frozen diffusion with constant immigration rate
$\sigma$, driven by the same white noise and with $Y_0^\sigma = x$.
By the comparison property in Proposition \ref{prop:logistic-comparison}, we have $Y_t^\sigma\leq X_t$ for every $t \ge 0$ almost surely. Recall that $\mbb{P}(\tau_0^{Y^\sigma} < \infty) = 0$ in this case. Then we have $\mathbb P(\tau_0<\infty)=0.$ Conversely, suppose that $am_{t_0} < \sigma$ for some $t_0 \ge 0$. By continuity, there exist $0 \le r < u$ and $0 \le \overline{c} < \sigma$ such that $am_t \le \overline{c}$ for $t \in [r, u]$. On $\{\tau_0 > r\}$, let $Y^{\overline{c}}$ be the frozen diffusion with immigration rate $\overline{c}$, started from $Y_r^{\overline{c}} = X_r$ and driven by the same white noise. Then we have $X_t \le Y^{\overline{c}}_t$ for every $t \in [r, u]$ almost surely. Conditioning on $\mathcal{F}_r$, and by the Markov property, we obtain
\beqnn
\mathbb{P}_x(\tau_0 \le u) \ge \mathbb{P}_x(\tau_0 \le r) + \mathbb{E}_x[1_{\{\tau_0 > r\}}p_{u-r}^{\overline{c}}(X_r)].
\eeqnn 
If $\mathbb{P}_x(\tau_0 \le r) > 0$, the claim is immediate. Otherwise, $\mathbb{P}_x(\tau_0 > r) = 1$ and $X_r > 0$ almost surely, so the expectation is strictly positive. Then (i) follows.

For (ii), assume that $\limsup_{t \rightarrow \infty}am_t < \sigma$. Then
there exist $T > 0$ and $c_1< \sigma$ such that $am_t \le c_1$ for all
$t \ge T$. In $\{\tau_0 > T\}$, conditionally on $\mathcal{F}_T$, let
$Y^{c_1}$ be the frozen diffusion with the immigration rate $c_1$ and the initial
value $Y_T^{c_1} = X_T$. Then we have $X_t \le Y_t^{c_1}$ for all
$t \ge T$ almost surely. Recall that $Y^{c_1}$ hits $0$ almost surely from
every positive initial value. We obtain 
$\mbb{P}_x(\tau_0 = \infty|\mathcal{F}_T) = 0$ on $\{\tau_0 > T\}$. Then
$\mbb{P}_x(\tau_0 = \infty) =
\mbb{E}_x[1_{\{\tau_0 > T\}}\mbb{P}_x(\tau_0 = \infty|\mathcal{F}_T)] = 0$.
When $a \le a_\star^{\rm log}$, one sees that $m_t \rightarrow 0$ as
$t \rightarrow \infty$ by Proposition \ref{prop:logistic-threshold}. Then
(ii) follows.

We next prove (iii). Recall that $F$, defined by \eqref{F}, is strictly
increasing from $a_\star^{\rm log}$ to infinity. Hence
$a_{\text{hit}} = \sigma I_1/I_2 = F(1) > a_\star^{\rm log}$. Moreover, when $a > a_\star^{\rm log}$,    and $F(\alpha_a) = a$,  then $\alpha_a < 1$ if and
only if $a < a_{\text{hit}}$. If
$a \le a_\star^{\rm log}$, we have $\mbb{P}_x(\tau_0 < \infty) = 1$ by (ii).
If $a_\star^{\rm log} < a < a_{\text{hit}}$, then $\alpha_a < 1$. By
Proposition \ref{p0814}, it yields that
$\lim_{t \rightarrow \infty}am_t = am_a = \sigma \alpha_a < \sigma$. Again
by (ii), we have $\mbb{P}_x(\tau_0 < \infty) = 1$. On the other hand, for
the case of $a > a_{\text{hit}}$, we have $\alpha_a > 1$. By Proposition
\ref{p0814}, there exists $T > 0$ such that $am_t \ge \sigma$ for every
$t \ge T$.
Let $Y^0$ be the frozen diffusion without immigration, starting from $Y^0_0 = x$ and driven by the same white noise as $X$. Then $\mbb{P}_x(\tau_0^{Y^0} < \infty) = 1$ and $\mbb{P}_x(\tau_0^{Y^0} \le T) \in (0, 1)$. Since $am_t \ge 0$, by Proposition \ref{prop:logistic-comparison}, we have  $X_t \ge Y^0_t$ for any $t \ge 0$ almost surely. Then $\mbb{P}_x(\tau_0 > T) \ge \mbb{P}_x(\tau_0^{Y^0} > T) > 0$.
 In event $\{\tau_0 > T\}$, conditionally on $\mathcal{F}_T$, let $Y^\sigma$ be the frozen diffusion with immigration rate $\sigma$ and $Y_T^\sigma = X_T$. Then $X_t \ge Y_t^\sigma$ for $t \ge T$ almost surely. As $0$ is inaccessible for $Y^\sigma$, $X$ cannot hit $0$ after $T$ on this event. Therefore, $\mathbb{P}_x(\tau_0 = \infty) \ge \mathbb{P}_x(\tau_0 > T) > 0$, and therefore $\mbb{P}_x(\tau_0 < \infty) < 1$. If, in addition, $am_t < \sigma$ for some $t \ge 0$, then we have $\mbb{P}_x(\tau_0 < \infty) > 0$. 
This completes the proof of (iii).
\qed

\section{Appendix}
\label{sec:wellposedness-proofs}

The following result is a time-inhomogeneous version of the comparison
principle of Dawson and Li~\cite[Theorem~2.2]{DawsonLi2012}. Since their
theorem is stated for time-homogeneous coefficients, we include the
short argument needed for locally bounded deterministic functions
$q_i$.

\begin{pro1}\label{prop:logistic-comparison}
Let $q_1,q_2:[0,\infty)\to[0,\infty)$ be locally bounded
Borel functions. For $i=1,2$, let $Z^i$ be a nonnegative
nonexplosive solution of
\beqnn
Z_t^i
=Z_0^i+\int_0^t
\bigl(q_i(s)-\beta Z_s^i-b_2(Z_s^i)^2\bigr)\,\d s
+\int_0^t\int_0^{Z_s^i}W(\d s,\d u),
\eeqnn
where the two equations are driven by the same Gaussian white noise
$W$ with intensity $2\sigma\,\d s\,\d u$. If $Z_0^1\leq Z_0^2$ almost surely and $q_1(t)\leq q_2(t)$ for any $t\geq0,$
then $\mathbb P\bigl(Z_t^1\leq Z_t^2
\text{ for every }t\geq0\bigr)=1.$
\end{pro1}

\proof
Set $D_t:=Z_t^1-Z_t^2$.
Then
\beqnn
D_t
=D_0+\int_0^t
\left[q_1(s)-q_2(s)-\beta D_s
-b_2D_s(Z_s^1+Z_s^2)\right]\d s+M_t,
\eeqnn
where 
the quadratic variation of $M$ satisfies
\beqnn
\d\langle M\rangle_t =
2\sigma\int_0^\infty\left|\mathbf 1_{\{u\leq Z_t^1\}}-\mathbf 1_{\{u\leq Z_t^2\}}\right|^2\,\d u\,\d t = 2\sigma|D_t|\,\d t.
\eeqnn
Consequently, the occupation-density formula and dominated convergence
give
\beqnn
L_t^0(D)
=
\lim_{\varepsilon\downarrow0}\frac{1}{2\varepsilon}
\int_0^t\mathbf 1_{\{|D_s|\leq\varepsilon\}}
\,\d\langle M\rangle_s =
\lim_{\varepsilon\downarrow0}\frac{\sigma}{\varepsilon}
\int_0^t\mathbf 1_{\{|D_s|\leq\varepsilon\}}
|D_s|\,\d s=0.
\eeqnn
For $n\geq1$, define $
\rho_n:=\inf\{t\geq0:Z_t^1+Z_t^2\geq n\}.$
Applying the It\^o--Tanaka formula
(see, e.g., \cite[Theorem~142]{Situ}) to $D^+$ at
$t\wedge\rho_n$, we obtain
\beqnn
D_{t\wedge\rho_n}^+
\ar=\ar D_0^+
+\int_0^{t\wedge\rho_n}\mathbf 1_{\{D_s>0\}}
\left[q_1(s)-q_2(s)-\beta D_s
-b_2D_s(Z_s^1+Z_s^2)\right]\d s\cr
\ar\ar
+\int_0^{t\wedge\rho_n}
\mathbf 1_{\{D_s>0\}}\,\d M_s,
\eeqnn
since $L_{t\wedge\rho_n}^0(D)=0$.  Let
$\beta^-:=(-\beta)\vee0$. Then on $\{D_s>0\}$, we have
\beqnn
q_1(s)-q_2(s)-\beta D_s
-b_2D_s(Z_s^1+Z_s^2)
\leq\beta^-D_s^+.
\eeqnn
Moreover, the stopped stochastic integral is a square-integrable
martingale because its quadratic variation is bounded by
$2\sigma nt$. Therefore, $\mathbb E[D_{t\wedge\rho_n}^+]
\leq\mathbb E[D_0^+]
+\beta^-\int_0^t
\mathbb E[D_{s\wedge\rho_n}^+]\,\d s.$
Since $D_0^+=0$ almost surely, by Gronwall's inequality, we have
$\mathbb E[D_{t\wedge\rho_n}^+]=0$. Nonexplosion implies that
$\rho_n\uparrow\infty$ almost surely and, hence, 
$\mathbb E[D_t^+]=0$ for every fixed $t\geq0$. Applying this conclusion
at all rational times and using path continuity gives
$D_t\leq0$ simultaneously for every $t\geq0$ almost surely.
\qed

{\it Proof of Proposition \ref{exist and unique}.}
 Let $T>0$ and $(m_t)_{0 \le t \le T} \in C([0,T],\R_+)$.  On the given stochastic basis, we consider the following frozen equation:
\beqlb\label{eq:frozen-mean-input}
X_t^m = X_0+\int_0^t\int_0^{X_{s}^m}W(\d s,\d u)
-\int_0^t[bX_s^m + g(X_s^m)]\,\d s +a\int_0^tm_s\,\d s
\eeqlb
for $t \in [0, T].$
By the localization construction in \cite[Theorem 2.5]{DawsonLi2012} and the Yamada-Watanabe principle (see, e.g., \cite[page 104]{Situ}), there exists a nonnegative strong solution $(X_t^m)_{0 \le t \le T}$ to \eqref{eq:frozen-mean-input}. The same conclusion holds when $(m_t)_{0 \le t \le T}$ is a random adapted process.
 Since $g\geq0$, it follows that
\begin{equation}\label{eq:frozen-first-moment}
\E [X_t^m]
\leq e^{b^-t}\E [X_0]
+a\int_0^te^{b^-(t-s)}m_s\,\d s, 
\qquad 0\leq t\leq T,
\end{equation}
where $b^- := (-b)\vee 0.$ Moreover, for $t \in [0, T]$,
\beqnn
\mathbb{E}\left[\int_0^t g(X_s^m) \d s\right] \le \mathbb{E}[X_0] + b^-\int_0^t \mathbb{E}[X_s^m]\d s + a\int_0^t m_s \d s,
\eeqnn 
which is finite by \eqref{eq:frozen-first-moment}.   
To compare two frozen equations, we choose numbers
$1=a_0>a_1>a_2>\cdots\downarrow0$ such that
$\int_{a_k}^{a_{k-1}}r^{-1}\,\d r=k$. Let
$\psi_k$ be  a continuous function supported in $(a_k,a_{k-1})$ such that $
0\leq\psi_k(r)\leq\frac{2}{kr}$ and $\int_{a_k}^{a_{k-1}}\psi_k(r)\,\d r=1,$
and let
\[
\phi_k(z)=\int_0^{|z|}\int_0^y\psi_k(r)\,\d r\,\d y.
\]
Then $0\leq\phi_k'(z)\sgn(z)\leq1$,
$|z|\phi_k''(z)\leq2/k$, and $\phi_k(z)\uparrow|z|$.
Let $X^m$ and $Y^n$ solve the frozen equations with $(m_t)_{0 \le t \le T}$ and $(n_t)_{0 \le t \le T}$, driven by the
same white noise, and with initial values $X_0$ and $Y_0$. Recall that $g$ is non-decreasing. Then $(g(x) - g(y))\sgn(x - y) \ge 0$. By It\^o's formula, one sees that
\begin{align*}
\mathbb{E}[\phi_k(X_t^m-Y_t^n)]
={}&\mathbb{E}[\phi_k(X_0-Y_0)] +\sigma\int_0^t\mathbb{E}\left[
\phi_k''(X_s^m-Y_s^n)|X_s^m-Y_s^n|\right]\d s\\
&+\int_0^t\mathbb{E}\left[
\phi_k'(X_s^m-Y_s^n)\{-b(X_s^m-Y_s^n)+a(m_s-n_s)\}\right]\d s\\
&-\int_0^t\mathbb{E}\left[
\phi_k'(X_s^m-Y_s^n)\{g(X_s^m)-g(Y_s^n)\}\right]\d s\\
\leq{}&\mathbb{E}[\phi_k(X_0-Y_0)]
+b^-\int_0^t\mathbb{E}[|X_s^m-Y_s^n|]\,\d s\\
&+a\int_0^t|m_s-n_s|\,\d s+\frac{2\sigma t}{k}.
\end{align*}
Letting $k\to\infty$ and using Gronwall's inequality, we obtain
\begin{equation}\label{eq:frozen-stability}
\E[|X_t^m-Y_t^n|]
\leq e^{b^-t}\E[|X_0-Y_0|]
+a\int_0^te^{b^-(t-s)}|m_s-n_s|\,\d s.
\end{equation}

We now consider the mean-field equation.  For $(m_t)_{0 \le t \le T}\in C([0, T], \mathbb{R}_+)$, we define $
(\Phi_Tm)_t:=\E f(X_t^m)$ for $t \in [0, T].$ 
By \eqref{bound of f} and \eqref{eq:frozen-first-moment}, for $0 \le s < t \le T$, we have
\beqnn
\mathbb{E}[|X_t^m - X_s^m|] \ar\le\ar \left(2\sigma\int_s^t \mathbb{E}[X_r^m] \d r\right)^{1/2} + |b|\int_s^t \mathbb{E}[X_r^m]\d r\cr
\ar\ar + \int_s^t\mathbb{E}[g(X_r^m)]\d r + a\int_s^t m_r \d r \rightarrow 0
\eeqnn
as $t - s \rightarrow 0$. Hence $X^m$ is $L^1$-continuous. By the Lipschitz property of $f$,  one sees that $\Phi_Tm \in C([0, T], \mathbb{R}_+)$. Thus, $\Phi_T$ maps
$C([0, T], \mathbb{R}_+)$ to itself. Let $\|h\|_\gamma:=\sup_{0\leq t\leq T}e^{-\gamma t}|h_t|$. For $\gamma>b^-$,   \eqref{eq:frozen-stability} gives 
\beqnn
\|\Phi_Tm-\Phi_Tn\|_\gamma \ar\le\ar \sup_{0 \le t \le T}e^{-\gamma t}\mathbb{E}[|f(X_t^m) - f(X_t^n)|] \le C_f \sup_{0 \le t \le T}e^{-\gamma t}\mathbb{E}[|X_t^m - X_t^n|]\cr
\ar\le\ar aC_f \sup_{0 \le t \le T}e^{-\gamma t}\int_0^t e^{b^-(t-s)}|m_s-n_s|\,\d s
\leq\frac{aC_f}{\gamma-b^-}\|m-n\|_\gamma.
\eeqnn
Choose $\gamma>b^-+aC_f$.  Then $\Phi_T$ is a contraction and has a unique fixed point.
The solution of the frozen equation for this fixed point satisfies \eqref{X}. The uniqueness of the fixed point makes these solutions compatible as $T$ varies, so the construction gives a global nonnegative strong solution.

If $X$ and $Y$ are two solutions of \eqref{X} with the same initial value and white
noise, apply \eqref{eq:frozen-stability} with
$m_t=\E [f(X_t)]$ and $n_t=\E [f(Y_t)]$.  The Lipschitz property of $f$ and Gronwall's
inequality give $\E[|X_t-Y_t|]=0$ for every $t$, and hence pathwise uniqueness.  For a weak
solution, the deterministic function $m_t=\E [f(X_t)]$ must be the unique fixed point of
$\Phi_T$. Weak uniqueness for the frozen equation therefore gives weak
uniqueness for \eqref{X}.

It remains to check the second moment.  The fixed point $m$ is bounded in every finite
time interval. Let $\tau_k := \inf\{t \ge 0: X_t \ge k\}.$ It should be noted that $\tau_k \rightarrow \infty$ almost surely as $k \rightarrow \infty$, which follows immediately from the continuity of the paths. By  It\^o's formula,
\beqnn
\mathbb{E}[X_{t\wedge\tau_k}^2] \ar\le\ar  \mathbb{E}[X_0^2] + \mathbb{E}\left[\int_0^{t\wedge\tau_k} 2X_s[-bX_{s} - g(X_s) + am_s + \sigma] \d s\right] \cr
\ar\leq\ar\mathbb{E} [X_0^2]+C_T\int_0^t\left[1+\mathbb{E}[X_{s\wedge\tau_k}^2]\right]\,\d s,
\qquad0\leq t\leq T.
\eeqnn
 Gronwall's inequality and Fatou's lemma yield, as $k \rightarrow \infty$,
$\sup_{0\leq t\leq T}\E [X_t^2]<\infty$. \qed

{\it Proof of Proposition~\ref{prop:particle-wellposedness}.}
 There exists a unique strong solution to the frozen equation \eqref{eq:frozen-mean-input} when $(m_t)_{0 \le t \le T}$ is a nonnegative
predictable locally integrable process, which follows from the localization argument,
since it is additive and does not enter the noise coefficient.

Let $f^N(x)=N^{-1}\sum_{j=1}^Nf(x^j)$.  We first construct the solution by Picard iteration. Let $X^{i,N,(0)}$ be the unique strong solution to \eqref{eq:frozen-mean-input} with white noise $W_i$ and $m_t = 0$. Now we define $X^{N,(k+1)}$ by
\beqnn
X_t^{i,N,(k+1)} \ar=\ar X_0^{i,N} +\int_0^t\int_0^{X_{s}^{i,N,(k+1)}}  W_i(\d s,\d u) -\int_0^t\{bX_s^{i,N,(k+1)} +g(X_s^{i,N,(k+1)})\}\d s\cr
\ar\ar +a\int_0^t f^N(X_s^{N,(k)})\d s, \qquad i = 1, \cdots, N.
\eeqnn
Recall that there exists a nonnegative strong solution to \eqref{eq:frozen-mean-input} when $(m_t)_{0 \le t \le T}$ is a random adapted process. For the above equation, one can obtain the existence and uniqueness of a nonnegative strong solution recursively for $k = 1, 2, \cdots$. Consequently, every iterate is well-defined. Define
\beqnn
U_{k+1}(t):=\frac1N\sum_{i=1}^N
\mathbb{E}\left[\left|X_t^{i,N,(k+1)}-X_t^{i,N,(k)}\right|\right].
\eeqnn
By \eqref{eq:frozen-stability} and the Lipschitz continuity of $f$, one sees that, for $k \ge 1$,
\beqnn
U_{k+1}(t) \leq aC_f\int_0^te^{b^-(t-s)}U_k(s) \d s.
\eeqnn
Iterating the preceding inequality yields $U_k(t)\leq A_T\frac{(c_Tt)^k}{k!}$ for $0 \le t \le T$, where  $c_T=aC_fe^{b^-T}$ and $A_T=\sup_{0\leq t\leq T}U_1(t)<\infty.$ One can check that the sequence $\{X^{N,(k)}: k \ge 1\}$ is Cauchy in $\frac1N\sum_{i=1}^N\mathbb{E}\sup_{0\leq t\leq T}|\cdot|.$ In fact,  the Burkholder-Davis-Gundy inequality gives
\beqnn
S_{k+1}(T) \ar:=\ar \frac1N \sum_{i = 1}^N \mathbb{E}\left[\sup_{0\le t \le T}|X_t^{i,N,(k+1)} - X_t^{i,N,(k)}|\right]\cr
\ar\le\ar C_T \left[\left(\int_0^T U_{k+1}(s)\d s\right)^{1/2} + \int_0^T [U_{k+1}(s) +U_k(s)]\d s\right],
\eeqnn
which implies that $\sum_{k \ge 1} S_{k+1}(T) < \infty$.
Let $X^N$ denote the continuous adapted limit of $X^{N,(k)}$. One can see that $X^N$ is a solution to \eqref{Xi}.

Let $V(x)=\sum_{i=1}^Nx_i$. By \eqref{generator of particle} and \eqref{bound of f} we have
\beqnn
\mathcal{L}_NV(x) \le aNf(0)+(b^-+aC_f)V(x)
\eeqnn
and
\beqnn
\mathcal{L}_N V^2(x) \ar=\ar 2V(x)\left[-bV(x)-\sum_{i=1}^Ng(x_i)+a\sum_{i=1}^Nf(x_i)\right]+2\sigma V(x)\cr
\ar\le\ar  2(b^-+aC_f)V^2(x)+2[aNf(0)+\sigma]V(x) \cr
\ar\le\ar C_N[1+V^2(x)].
\eeqnn
Notice that $V(X_t^N) = \sum_{i = 1}^N X_t^{i,N}$.
Let $\tau_k^N :=\inf\{t \ge 0: V(X_t^N) \ge k\}$. By It\^{o}'s formula, we have
\beqnn
\mathbb{E}[V^2(X_{t\wedge\tau_k^N}^N)] \ar=\ar  \mathbb{E}[V^2(X_0^N)] +  \mathbb{E}\left[\int_0^{t\wedge\tau_k^N} \mathcal{L}_N V^2(X_s^N)\d s\right]\cr
\ar\le\ar  \mathbb{E}[V^2(X_0^N)] + C_N\int_0^t [1 + \mathbb{E}[V^2(X_{s\wedge\tau_k^N}^N)]]\d s.
\eeqnn
By Gronwall's inequality and letting $k \rightarrow \infty$, we have
\beqnn
\sup_{0 \le t \le T}\mathbb{E}[V^2(X_t^N)] = \sup_{0 \le t \le T}\mathbb{E}\left[\left(\sum_{i = 1}^N X_t^{i,N}\right)^2\right] < \infty.
\eeqnn

It remains to prove pathwise uniqueness.  Let $X^N$ and $Y^N$ be two solutions driven by the same initial value and white noise. For $\phi_k$ defined in the proof of Proposition \ref{exist and unique}, by It\^{o}'s formula and the Lipschitz continuity of $f$,
\beqnn
\mathbb{E}[\phi_k(X_t^{i,N} - Y_t^{i,N})] \ar\le\ar b^-\int_0^t \mathbb{E}[|X_s^{i,N} - Y_s^{i,N}|] \d s\cr
\ar\ar + a\int_0^t \mathbb{E}\left[\left|f^N(X_s^N) - f^N(Y_s^N)\right|\right]\d s + \frac{2\sigma t}{k}\cr
\ar\le\ar b^-\int_0^t \mathbb{E}[|X_s^{i,N} - Y_s^{i,N}|] \d s\cr
\ar\ar + aC_f\int_0^t \frac{1}{N}\sum_{j = 1}^N\mathbb{E}[|X_s^{j,N} - Y_s^{j,N}|]\d s + \frac{2\sigma t}{k}.
\eeqnn
Letting $k \rightarrow \infty$ and summing up from $i$ to $N$, it follows that
\beqnn
\frac{1}{N}\sum_{i=1}^N\E\left[|X_t^{i,N}-Y_t^{i,N}|\right] \le (b^- + aC_f) \int_0^t \frac{1}{N}\sum_{i=1}^N\E\left[|X_s^{i,N}-Y_s^{i,N}|\right] \d s,
\eeqnn
which implies that $\frac{1}{N}\sum_{i=1}^N\E\left[|X_t^{i,N}-Y_t^{i,N}|\right] = 0$ by  Gronwall's inequality.  The pathwise uniqueness follows from the continuity. \qed

{\bf Acknowledgements.}
The research of Shukai Chen is supported
    by the National Key R\&D Program of China
(No.2022YFA1006003), NSFC grant of China
(No.12401167),
Fujian Provincial Natural Science Foundation of China (No.2024J08050) and the Education and Scientific Research Project for Young and Middle-aged Teachers in Fujian Province of China
(No.JAT231015). The research of Lina Ji is supported in part by NSFC grant (No. 12301167), Guangdong Basic and Applied Basic Research Foundation (No. 2022A1515110986) and Shenzhen National Science Foundation (No. 20231128093607001). Xiaowen Zhou's research is supported
by the Natural Sciences and Engineering Research Council of Canada (RGPIN-2026-07619).

	\vskip 0.2truein
	\vskip 0.2truein

\bigskip

\noindent{\bf Shukai Chen:}  School of Mathematics and Statistics, Fujian Normal University,
Fuzhou 350117, Fujian, P. R. China. Email: {\texttt skchen@fjnu.edu.cn}

\bigskip

\noindent{\bf Lina Ji:}
MSU-BIT-SMBU Joint Research Center of Applied Mathematics, Shenzhen MSU-BIT University,
Shenzhen 518172, P. R. China. Email: {\texttt jiln@smbu.edu.cn}

\bigskip

\noindent{\bf Xiaowen Zhou:}
Department of Mathematics and Statistics, Concordia
			University, Montreal, Canada. Email: {\texttt xiaowen.zhou@concordia.ca.}

\end{document}